\documentstyle[12pt]{article}
  \newcommand{\const}{\rm const}
  \newcommand{\Var}{\rm Var}
  \newcommand{\supp}{\rm supp}
  \newcommand{\Law}{\rm Law}

  \newcommand{\Dom}{\rm  Dom}
  
  \newcommand{\card}{\rm card}
  \newcommand{\argmin}{\rm argmin}

\begin{document}

   \begin{center}

  \ \  \ \   {\bf  Moment and tail estimates and Banach space valued }\\

\vspace{4mm}

  \ \ \ \ \  {\bf Non-Central Limit Theorem (NCLT) for sums of  } \\

\vspace{4mm}

\ \ \ \  \  \ \ \ \ {\bf multi-indexed random variables, processes and fields }\\

\vspace{4mm}

  \ \ \ \  \ \ \ \ \ \  {\bf Ostrovsky E., Sirota L.}\\

\vspace{4mm}

  Bar-Ilan University, department of Mathematic and Statistics, ISRAEL, 59200. \\

\vspace{4mm}

E-mails: eugostrovsky@list.ru, \hspace{5mm} sirota3@bezeqint.net \\

\vspace{5mm}

  {\bf Abstract} \\

\vspace{4mm}

 \end{center}

\ We derive in this  preprint  the  moment and exponential tail estimates,
sufficient conditions for the Non-Central Limit Theorem (NCLT) in the ordinary
one-dimensional space as well as in the space of  continuous
functions for the  properly (natural) normalized multi-indexed sums of function of random  variables,  processes or fields (r.f.),
on the other words $ \ V \ - \ $ statistics, parametric, in general case. \par

 \ We construct also some examples in order to show the exactness of obtained estimates.\par

 \ We will use  the theory of the so-called  degenerate approximation of the functions of several variables as well as the theory of
 Grand Lebesgue Spaces (GLS) of measurable  functions  (random variables). \par

\vspace{5mm}

{\it Key words and phrases:}  Measure and probability, measurable functions, random variable and vector (r.v.), normalized sum, moment
and tail estimates, Non-Central Limit Theorem (NCLT), Banach space of continuous functions, degenerate functions and
approximation,  triangle inequality,  tail of distribution, Lebesgue-Riesz, Orlicz and Grand Lebesgue Spaces (GLS), distance function,
compact metric space, metric entropy and entropy integral, cardinal number of the set,
Dharmadhikari-Jogdeo-Rosenthal's constant and inequality, generating function, Lyapunov's inequality, Young-Orlicz function,
conditional expectation, Young-Fenchel transform, white Gaussian measure, rearrangement invariant (r.i.) space. \par

\vspace{4mm}


\vspace{5mm}

\section{ Definitions. Notations. Previous results. Statement of problem.}

\vspace{4mm}

 \ Let $ \   (X, \cal{B},  {\bf \mu})   \ $ and $ \   (Y, \cal{C}, {\bf \nu} ) \  $
be  two probability spaces:  $ \ \mu(X)  = \nu(Y) = 1 . \ $ We will denote by $ \  |g|_p  = |g|L(p) \ $ the ordinary
Lebesgue-Riesz $ \  L(p) \  $ norm of arbitrary measurable numerical valued function $ \   g: X \to R: \  $

$$
|g|_p = |g| L(p) = |g| L_p(X, \mu) := \left[ \int_X |g(x)|^p \ \mu(dx)  \right]^{1/p}, \  p \in [1, \infty)
$$
 analogously for the (measurable) function $ \  h: Y \to R \ $

$$
|h|_p = |h|L(p) = |h|L_p(Y,\nu) := \left[ \int_Y |h(y)|^p \ \nu(dy) \right]^{1/p};
$$
and for arbitrary integrable function of two variables $ \  f: X \otimes Y  \to R \ $

$$
|f|_p  = |f|L(p) = |f|L_p(X,Y):= \left[ \int_X \int_Y |f(x,y)|^p \ \mu(dx) \ \nu(dy) \right]^{1/p}, \ p \in [1, \infty).
$$

\vspace{4mm}

 \ Let $ \  Z_+ = \{  1,2, 3, \ldots \} \  $ and denote $ \ Z^2_+ = Z_+ \otimes Z_+, \  Z^d_+ = \otimes_{k=1}^d Z_+. \ $
Let also  $ \  \{\xi(i) \} \  $ and $ \  \{\eta(j)\},
 \ i,j = 1,2,\ldots, \   \xi := \xi(1), \ \eta := \eta(1)  $ be {\it common} independent random variables  defined on certain probability space
$ \   (\Omega, \cal{M}, {\bf P})  \ $  with distributions  correspondingly $ \ \mu, \ \nu: \ $

$$
{\bf P}(\xi(i)  \in A) = \mu(A), \ A \in \cal{B};
$$

$$
{\bf P}(\eta(j) \in F) = \nu(F), \ F \in \cal{C}, \eqno(1.0)
$$

so that

$$
{\bf E} |g(\xi)|^p  = |g|_p^p, \ {\bf E} |h(\eta)|_p^p = |h|^p_p
$$
and

$$
{\bf E} |f(\xi, \eta)|^p = |f|^p_p.
$$

 \ Let also $ \  L \ $ be arbitrary non - empty {\it finite} subset of the set $ \  Z^2_+; \ $ denote by $  \   |L| \  $  a  numbers of its
elements (cardinal number): $  \  |L| := \card(L). \  $  It is reasonable  to suppose in what follows $ \ |L| \ge 1. \ $ \par

 \ Define for any {\it  centered} function $ \  f: X \otimes Y \to R, \  $ i.e. for which

$$
  {\bf E}f(\xi, \eta) = \int_X \int_Y f(x,y) \ \mu(dx) \ \nu(dy) = 0,
$$
the following normalized sum

$$
S_L[f] \stackrel{def}{=} |L|^{-1/2} \sum_{(k(1), k(2)) \in L} f(\xi(k(1), \eta(k(2)), \eqno(1.1)
$$
which is a slight generalization of the classical $ \ U \ $ and $ \  V \ $ statistics, see the classical  monograph of
Korolyuk V.S and Borovskikh Yu.V. [24].  Offered here report is the direct generalization of a recent article [33], but we
apply here other methods.\par

 \ The {\it  reasonableness } of this norming function $ \   |L|^{-1/2} \ $
implies that in general, i.e. non-degenerate case $ \   \Var(S_L) \asymp 1, \ |L| \ge 1.  \ $ This propositions
holds true still in the multidimensional case. \par

 \ Our notations  and some previous results are borrowed from the works of S.Klesov  [20]-[24]. \par

\vspace{4mm}

 \ {\bf   Our claim in this report is to derive  the moment and exponential bounds for tail of distribution for the normalized sums of
multi-indexed independent random variables from (1.1). } \par
 {\bf We deduce also as a consequence the sufficient conditions for the weak compactness in the space of continuous function the
introduces sums in the case when the function $ \ f \ $ dependent still on some parameter  (parameters).   } \par

\vspace{4mm}

 \ Offered here results are generalizations of many ones obtained by S.Klesov  in an articles [20]-[24],  see also  the articles of N.Jenish,
and I.R.Prucha  [19], and  M.J.Wichura  [44], where was obtained in particular the CLT for these sums. \par

\vspace{4mm}

 \ The multidimensional  case, i.e. when $ \ \vec{k} \in  Z_+^d, \ $ will be considered further. \par

\vspace{4mm}

 \ The paper is organized as follows. In the  second section we describe and investigate the notion of
 the degenerate functions and approximation. In the next section we obtain one of the main results:
 the  moment estimates for multi-index sums. We outline in the fourth section the so-called
 Non-Central Limit Theorem for the multi-indexed sums. \par
 \  The fifth section contains the  multidimensional generalization of obtained results. The sketch of the theory of
Grand Lebesgue Spaces (GLS) with some new facts is represented in the next section. \par
 \ The exponential estimates for tail of multi-index sums is the content of $ \ 7^{th} \ $ section.The $ \ 8^{th} \ $ one
is devoted to the so-called  Banach space valued Non-Central Limit Theorem for multi-index sums.   \par
 \ The last section contains as ordinary some concluding remarks. \par

\vspace{4mm}

 \section{Degenerate functions and approximation. }

\vspace{4mm}

 \ {\bf Definition 2.1.}  The measurable centered function $ \  f:  X \otimes Y \to R  \ $ is said to be {\it  degenerate, } if it has a form

$$
f(x,y) = \sum_{i = 1}^M \sum_{j=1}^M \lambda_{i,j} \ g_i(x) \ h_j(y),  \eqno(2.1)
$$
where $ \  \lambda_{i,j} = \const \in R, \  M =  \const = 1,2, \ldots, \infty. \ $ \par

 \ The degenerate functions (and kernels) of the form (2.1) are used, e.g., in the approximation theory, in the theory of random processes and fields,
in the theory of integral equations, in the game theory etc. \par

 \  A particular application of this notion may be found in the  authors article [37]. \par

\vspace{4mm}

 \ {\it  It will be presumed in this report in addition to the expression (2.1) that all the functions   \ } $ \  \{  g_i \}, \ \{ h_j \}  \  $
{\it are centered, such that correspondingly}

$$
{\bf E } g_i(\xi) = \int_X g_i(x) \ \mu(dx) = 0 \eqno(2.2a)
$$
and
$$
{\bf E} h_j(\eta) = \int_Y h_j(y) \ \nu(dy) =  0; \ i,j = 1,2,\ldots,M. \eqno(2.2b)
$$

\vspace{4mm}

 \ Denotation: $ \ M = M[f] \stackrel{def}{=}  \deg(f); \ $ of course, as a capacity of the value $ \ M \ $  one can understood its {\it constant}
minimal value. \par

 \ Two examples. The equality (2.1) holds true if the function $ \  f(\cdot, \cdot) \  $  is trigonometrical or algebraical polynomial.  \par
\ More complicated example: let $  \   X \ $ be compact metrizable space  equipped with the non-trivial probability Borelian measure $ \ \mu. \ $
This  imply that an arbitrary non-empty open set has a positive measure. \par

 \ Let also $ \  f(x,y), \  x,y \in X \ $ be continuous numerical valued  non-negative definite function. One can write the famous Karunen-Loev's
decomposition

$$
f(x,y) = \sum_{k=1}^M  \lambda_k  \ \phi_k(x) \ \phi_k(y),
$$
where $ \  \lambda_k, \ \phi_k(x) \ $ are  correspondingly eigenvalues and  orthonormal eigenfunction  for the function (kernel)  $ \ f(\cdot, \cdot): \ $

$$
\lambda_k \ \phi_k(x)  = \int_X f(x,y) \phi_k(y) \ \mu(dy). \eqno(2.3)
$$
 \ We assume without loss of generality

$$
\lambda_1 \ge \lambda_2 \ge \ldots \lambda_k \ge \ldots \ge 0.  \eqno(2.3a)
$$

\vspace{4mm}

\  Further, let  $ \  B_1, \ B_2, \  B_3, \ \ldots, B_M \  $  be some rearrangement invariant  (r.i.) spaces builded correspondingly over the
spaces $ \ X,Y; \ Z,W, \ldots, \ $
for instance, $ \ B_1 = L_p(X), \  B_2 = L_q(Y), 1 \le p,q \le \infty.  \  $ If $ \  f(\cdot) \in B_1 \otimes B_2, \  $ we suppose also
in (2.1) $ \ g_i \in B_1, \ h_j \in B_2; \ $ and if in addition in (2.1) $ \  M = \infty,  \  $ we suppose that the series in (2.1)
converges in the norm $  \  B_1 \otimes B_2 \   $

$$
\lim_{m \to \infty} || \   f(\cdot) -   \sum_{i = 1}^m \sum_{j=1}^m \lambda_{i,j} \ g_i(\cdot) \ h_j(\cdot)  \  ||B_1 \otimes B_2 = 0.
\eqno(2.3b)
$$

\  The condition (2.3b) is satisfied if for example $ \ ||g_i||B_1  = ||h_j||B_2 = 1 \ $ and

$$
 \sum \sum_{i,j = 1}^M |\lambda_{i,j}| < \infty, \eqno(2.4)
$$
or more generally when

$$
\sum \sum_{i,j = 1}^M |\lambda_{i,j}| \cdot ||g_i||B_1 \cdot ||h_j||B_2 < \infty. \eqno(2.4a)
$$

 \ The function of the form (2.1) with $ \    M = M[f] = \deg (f) < \infty \ $ is named {\it degenerate},
notation $ \  f \in D[M]; \ $ we put also $  \  D := \cup_{M < \infty} D[M]. \ $ Obviously,

$$
B_1 \otimes B_2 = D[\infty].
$$

 \ Define also for each such a function $ \   f \in D \  $ the following non-negative quasi-norm

$$
||f|| D(B_1, B_2) \stackrel{def}{=} \inf \left\{ \  \sum \sum_{i,j = 1,2,\ldots,M[f] }  |\lambda_{i,j}| \ ||g_i||B_1 \ ||h_j||B_2 \  \right\}, \eqno(2.5)
$$
where all the arrays $ \ \{ \lambda_{i,j} \} , \ \{ g_i\}, \ \{h_j\}  \ $ are taking from the  representation  2.1. \par

 \ We will write for brevity $ \ ||f||D_p   := $

$$
||f|| D(L_p(X), L_p(Y)) =
\inf \left\{ \  \sum \sum_{i,j = 1,2,\ldots,M[f] }  |\lambda_{i,j}| \ |g_i|_p \ |h_j|_p \  \right\}, \eqno(2.5a)
$$
where all the arrays $ \ \{ \lambda_{i,j} \} , \ \{ g_i\}, \ \{h_j\}  \ $ are taking from the  representation  2.1. \par

\vspace{4mm}

\ Further, let the function  $ \  f \in B_1 \otimes B_2  \ $ be given. The error of a degenerate approximation of the function $ \ f \ $
by the degenerate ones of the degree $ \ M \ $ will be introduced as follows

$$
Q_M[f](B_1 \otimes B_2) \stackrel{def}{=} \inf_{\tilde{f} \in D[M]} ||f - \tilde{f}||B_1 \otimes B_2 =
\min_{\tilde{f} \in D[M]} ||f - \tilde{f}||B_1 \otimes B_2. \eqno(2.6)
$$
 \  Obviously, $ \  \lim Q_M[f] (B_1 \otimes B_2) = 0, \ M \to \infty. \  $\par

 \ For brevity:

$$
Q_M[f]_p \stackrel{def}{=} Q_M[f](L_p(X) \otimes L_p(Y)). \eqno(2.6a)
$$

\vspace{4mm}

 \  The function $ \ \tilde{f} \ $ which realized the minimum in (2.6),  not necessary to be unique,
will be  denoted by $ \  Z_M[f](B_1 \otimes B_2):  \  $

$$
 Z_M[f](B_1 \otimes B_2):= \argmin_{\tilde{f} \in D[M]} ||f - \tilde{f}||B_1 \otimes B_2, \eqno(2.7)
$$
so that

$$
Q_M[f](B_1 \otimes B_2) = ||f - Z_M[f]||B_1 \otimes B_2.  \eqno(2.8)
$$

 \ For brevity:

$$
Z_M[f]_p := Z_M [f](L_p(X) \otimes L_p(Y)). \eqno(2.9)
$$

 \ Let for instance again $ \  f(x,y), \  x,y \in X \ $ be continuous numerical valued  non-negative definite function,
see (2.3) and (2.3a). It is easily to calculate

$$
Q_M[f] (L_2(X) \otimes L_2(X)) = \sum_{k=M + 1}^{\infty} \lambda_k.
$$

\vspace{4mm}

 \section{Moment estimates for multi-index sums.}

\vspace{4mm}

{\bf 0. \  Trivial estimate.  }  \par

\vspace{4mm}

 \ The following simple estimate based only on the triangle inequality, may be interpreted as trivial:

$$
|S_L|(L_p(X) \otimes L_p(Y)) \le |L|^{1/2} \ |f|L_p(X) \otimes L_p(Y)). \eqno (3.0)
$$

\vspace{4mm}

{\it  Hereafter \   } $ \ p \ge 2. \ $ \par

\vspace{4mm}

 \ {\bf 1. }  Let us consider at first the one-dimensional case, namely $ \ f = f(x) = g(x), \ x \in X: $

\vspace{4mm}

$$
S_n = S_n[f] := n^{-1/2} \sum_{i=1}^n g(\xi_i), \eqno(3.1)
$$
where as before $  \ {\bf E} g(\xi) = 0. \ $ One can apply the famous Dharmadhikari-Jogdeo-Rosenthal's inequality

$$
\sup_n  | S_n|_p \le K_R(p) \cdot |g(\xi)|_p, \ p \ge 2, \eqno(3.2)
$$
where the optimal {\it Rosenthal's  function} on $ \ p: \ p \to  K_R(p) \ $  in (3.2) may be estimated as follows: $ \  K_R(2) = 1 \  $ and

$$
K_R(p) \le C_R \ \frac{p}{ e \cdot \ln p} =: K_R^o(p),
$$
and  the exact value of {\it Rosenthal's constant} $ \ C_R \ $ is following

$$
C_R \approx 1.77638\ldots
$$
and this value is attained when $ \  p = p_0 \approx 33.4610\ldots,  \  $  see [17], [32].  \par

 \ Analogously

$$
\sup_n \left| n^{-1/2} \sum_{i=1}^n g(\xi_i) \  h(\eta_i)  \right|_p  \le K_R(p) \ |g(\xi)|_p \ |h(\eta)|_p, \ p \ge 2.
$$

\vspace{4mm}

{\bf  2.  The two-dimensional case.  } \par

\vspace{4mm}

 \ In this subsection  the kernel-function $ \  f = f(x,y) \  $ will be presumed to be degenerate with minimal degree $ \ M = M[f] = 1: \ $

$$
f(x,y) = g(x) \cdot h(y), \ x \in X, \ y \in Y.
$$

 \ Let us consider the correspondent double sum $  \  S_L[f] =  S^{(2)}_L :=   $

$$
 |L|^{-1/2} \sum  \sum_{i,j \in L}  g(\xi_i) \ h(\eta_j), \
 n =  \vec{n} = (n_1, n_2) \in L, \ n_1, n_2 \ge 1, \eqno(3.3)
$$
 where as before $ \  L \ $ is arbitrary non - empty: $ \ |L| \ge 1 \ $ subset of the integer positive plane $ \  Z^2_+.  \  $
O.Klesov in [22], [23] proved the following estimate

$$
\sup_{L: |L| \ge 1} \left| S^{(2)}_L   \right|_p \le K_R^2(p)  \ |g(\xi)|_p \ |h(\eta)|_p, \ p \ge 2, \eqno(3.4)
$$
and analogously for the multi-index sums. \par

\vspace{4mm}

{\bf  3.  Estimation for arbitrary degenerate kernel.  } \par

\vspace{4mm}

 \ In this subsection $ \ f(\cdot,\cdot) \ $ is degenerate:

$$
f(x,y) = \sum \sum_{k_1, k_2 = 1}^M \ \lambda_{k_1, k_2} \ g_{k_1}(x) \ h_{k_2}(y), \eqno(3.5)
$$
where

$$
g_{k_1}(\cdot) \in L_p(X), \ h_{k_2} (\cdot) \in L_p(Y),
$$
and as before $ \  {\bf E} g_{k_1}(\xi) = {\bf E} h_{k_2}(\eta) = 0.  \  $

 \ Let us investigate the introduced before  statistics

$$
S^{(\lambda)}_{L} = S^{(\lambda)}_L[f]  := | L|^{-1/2} \ \sum \sum_{i,j \in L} f(\xi_i, \eta_j) =
$$

$$
| L|^{-1/2} \ \sum \sum_{i,j \in L}  \ \sum \sum_{k_1, k_2 = 1}^M \lambda_{k_1, k_2} g_{k_1}(\xi_i) \ h_{k_2}(\eta_j) =
$$

$$
|L|^{-1/2} \ \sum \sum_{k_1, k_2 = 1}^M \lambda_{k_1, k_2} \ \sum \sum_{i,j  \in L} g_{k_1}(\xi_i) \ h_{k_2}(\eta_j), \ M < \infty. \eqno(3.6)
$$

 \  We have using the triangle inequality and the estimate (3.4)

$$
 \sup_{ L: |L| \ge 1} \left| S^{(\lambda)}_L[f] \right|_p \le  K_R^2(p) \cdot
\sum \sum_{k_1, k_2 = 1}^M \ \left| \lambda_{k_1, k_2} \right| \ |g_{k_1}(\xi)|_p  \ |h_{k_2}(\eta)|_p. \eqno(3.7)
$$

 \ This estimate remains true in the case when $ \ M = \infty, \ $ if of course the right-hand side of (3.7) is finite.  Therefore,

$$
 \sup_{ L: |L| \ge 1} \left| S^{(\lambda)}_L[f] \right|_p \le   K_R^2(p) \cdot D_p[f], \eqno(3.8)
$$
see (2.5), (2.5a). \par

\vspace{4mm}

{\bf  4. Main result. Degenerate approximation approach. } \par

\vspace{4mm}

 \ {\bf Theorem 3.1.}  Let $ \  f = f(x,y)  \  $ be arbitrary function from the space $ \   L_p(X) \otimes L_p(Y), \ p \ge 2. \  $ Then
$ \ \sup_{L: |L| \ge 1} |S_L[f]|_p \le W[f](p), $ where $ \ W[f](p) \stackrel{def}{=} $

$$
 \sup_{L: |L| \ge 1} \inf_{M \ge 1} \left[ \ K^2_R(p) \ || Z_M[f] ||D_p  + |L|^{1/2} Q_M[f]_p \  \right]. \eqno(3.9)
$$

\vspace{4mm}

{\bf Proof } is very simple, on the basis of previous results of this section. Namely, let $ \ L \ $ be an arbitrary non-empty set. Consider a
splitting

$$
f = Z_M[f] + ( f - Z_M[f] ) =: \Sigma_1 + \Sigma_2.
$$
 \ We have

$$
|\Sigma_1|_p =  |Z_M[S_L[f]] |L_p(X \otimes Y)   \le K_R^2(p) \ ||Z_M[f]||D_p.
$$
The member $  \ |\Sigma_2|_p $ may be estimated by virtue of inequality (3.0):

$$
|\Sigma_2|_p \le |L|^{1/2} \ |f - Z_M[f]|_p = |L|^{1/2} Q_M[f]_p.
$$
 \ It remains to apply the triangle inequality and minimization over $ \ M. \ $ \par

\vspace{4mm}

 \ {\bf Example 3.1.}  We deduce from (3.9) as a particular case

$$
 \sup_{ L: |L| \ge 1} \left| S^{(\lambda)}_L[f] \right|_p \le   K_R^2(p) \cdot ||f||D_p, \eqno(3.10)
$$
 if of course the right-hand side  of (3.10) is finite for some value $ \ p, \ p \ge 2. \ $ \par
\ Recall that in this section $ \ d = 2. \ $ \par

\vspace{4mm}

 \section{Non-Central Limit Theorem for multi-indexed sums.}

\vspace{4mm}

 \ We intend to find in this section the limit  in distribution as $  \ |L| \to \infty \ $ of the random value $ \  S_L[f].  \ $
The limit distribution coincides with a multiple stochastic integrals relative a white Gaussian measure, alike ones for $ \  U \ $ and
$ \ V \ $ statistics,  but applied in this preprint methods are essentially other, namely they based on the degenerate approximation
 of the kernel $ \ f(\cdot, \cdot). $ \par
 \ {\it  Additional assumptions and notations.  \ }  We suppose for beginning that the domain $ \ L \ $ is {\it rectangle: }

$$
L = L(n_1, n_2) = \{ (i,j): 1 \le i \le n_1; 1 \le j \le n_2  \}, \eqno(4.1)
$$
where both the boundaries $ \  n_1, \ n_2  \ $ tend to infinity: $ \ \lim n_s =\infty, \ s = 1,2. \ $ \par
 \ Further, we suppose   $ \ 0 < W[f](2) <  \infty; \ $ then the  family of distributions $ \  \Law(S_L) \ $ is weakly compact on the real line. \par

 \ Let  again the function $ \  f = f(x,y) \  $ be degenerate:

$$
f(x,y) = \sum \sum_{k_1, k_2 = 1}^M \ \lambda_{k_1, k_2} \ g_{k_1}(x) \ h_{k_2}(y), \ M \le \infty;  \eqno(4.2)
$$
but we impose on the system of a functions $ \  \{g_{k_1}, \ h_{k_2} \} \ $ without  loss of generality the following condition
of orthonormality:

$$
{\bf E} g_{k}(\xi) \ g_{l}(\xi) = {\bf E} h_{k}(\eta) \ h_{l}(\eta) = \delta_{k}^{l}, \eqno(4.3)
$$
where $ \ \delta_{k}^{l} \ $ is the Kronecker's symbol. \par

 \ We deduce from the condition $ \ 0 < W[f](2) <  \infty \ $ that

$$
\sigma^2 \stackrel{def}{=}  \Var [ f(\xi, \eta) ]  = \sum \sum_{k_1, k_2 = 1}^M \lambda^2_{k_1, k_2} \in (0, \infty). \eqno(4.4)
$$

 \ Introduce also two independent series of independent standard normal distributed random variables

$$
\Law(\tau_k) = N(0,1), \  \Law(\theta_k) = N(0,1).
$$

\vspace{4mm}

{\bf Theorem 4.1.} We assert under formulated above in this section  conditions that the sequence of the r.v. $ \  S_L \ $
converges in distribution as $ \  \min(n_1, \ n_2) \to \infty  \ $ to the multiple stochastic integral  $ \ S_{\infty} \ $ relative the
white Gaussian random measure with non-random square integrable integrand:

$$
S_L \stackrel{d}{\to} S_{\infty} \stackrel{def}{=} \sum \sum_{k,l = 1}^M \lambda_{k,l} \  \tau_k  \ \theta_l, \eqno(4.5)
$$
where the symbol $ \  S_L \stackrel{d}{\to} \nu \ $ denotes the weak (i.e. in distribution) convergence. \par

\vspace{4mm}

 \ {\bf Proof.} It is sufficient taking into account the weak compactness of $ \ S_L \ $
to consider the case when $ \ M = 1, \ $ as long as the general case is linear combination of the
one-dimensional ones. So, let

$$
f(x,y) = g(x) \cdot  h(y),
$$
where the functions $ \  g(\xi), h(\eta)  \ $ are centered and orthonormal. We  have as $ \  n_1, \ n_2 \to \infty \ $

$$
S_L = (n_1, n_2)^{-1/2} \sum_{i=1}^{n_1} \sum_{j=1}^{n_2} g(\xi_i) \ h(\eta_j) =
$$

$$
\left[ n_1^{-1/2} \sum_{i=1}^{n_1} g(\xi_i)  \right] \cdot \left[ n_2^{-1/2} \sum_{j=1}^{n_2} h(\eta_j)   \right]
\stackrel{d}{\to} \tau_1 \cdot \theta_1,
$$
on the basis of the classical one-dimensional CLT. Therefore, $ \  S_L \stackrel{d}{\to} S_{\infty}. \  $ \par
 \ In the general case:

$$
S_L \stackrel{d}{\to} S_{\infty} = \sum \sum_{k,l = 1}^M \lambda_{k,l} \  \tau_k  \ \theta_l, \eqno(4.6)
$$
since

$$
\sum \sum_{k,l = 1}^M \lambda_{k,l}^2 < \infty.
$$

 \ Notice that the right-hand side of (4.6) is {\it homogeneous} polynomial of degree 2 on the centered Gaussian
distributed independent r.v. Let us introduce thence the following white Gaussian measure $  \  \zeta \  $ with independent
values on the disjoint sets defined on the Borelian sets of positive plane $  \  R_+^2: \  $

$$
{\bf E} \ \zeta(A_1, B_1) \ \zeta(A_2, B_2) = \int_{A_1 \cap A_2} dx \ \int_{B_1 \cap B_2} dy.
$$

 \ All the r.v. $ \  \{  \tau_k, \ \theta_l  \} \  $ may be realized as follows.

$$
\tau_k := \int_{[k, k+1)} \int_0^1 \zeta(dx, dy); \ \theta_l := \int_2^3 \int_{[l, l+1)} \zeta(dx, dy);
$$
then the limiting r.v. $ \ S_{\infty} \ $ may be represented in the form

$$
S_{\infty}  = \sum \sum_{k,l = 1}^M \lambda_{k,l} \  \tau_k  \ \theta_l =  \sum \sum_{k,l = 1}^M \lambda_{k,l} \times
$$

$$
\int \int_{R^2}  \left[ I(x \in [k, k+1)) \ I(y \in [0,1]) + I(x \in [2,3]) \ I(y \in [l, l+1))    \right] \ \zeta(dx, dy) =
$$

$$
 \int \int_{R^2} r(x,y) \ \zeta(dx, dy), \eqno(4.7)
$$
where the square integrable on the whole plane deterministic (non-random) function $ \  r = r(x,y)  \  $ has a
form $ \  r(x,y) = \sum \sum_{k,l = 1}^M \ \lambda_{k,l} \times  $

$$
 \left[ \ I(x \in [k, k+1)) \ I(y \in [0,1]) +
 I(x \in [2,3]) \ I(y \in [l, l+1))    \right],\eqno(4.8)
$$
where $ \  I(z \in C) \  $ denotes the ordinary indicator function of the set $ \ C. \ $ \par
 \ So, the limit distribution $ \ \Law(S_L) \ $ of the sequence $ \  S_L \ $ coincides with the multiple stochastic integral (4.7) relative the
 centered Gaussian random  white measure. \par

 \vspace{4mm}

 \ {\it Let us consider a more general case, i.e. when the sets $ \ L \ $ are not rectangle.   } \par

 \vspace{4mm}

 \ We need to introduce some new geometrical notations. Denote by $ \  \pi_-(L) \  $ the set of all rectangles which are {\it inscribed}
into the set $ \ L: \ \pi_-(L) = \{ L_-  \}, $  where

$$
L_- = \{ [n(1)_{--}, n(1)_-] \otimes [n(2)_{--}, n(2)_- ]  \}: \ L_- \subset L, \eqno(4.9)
$$

$$
1 \le n(1)_{--}  \le n(1)_- < \infty, \ n(1)_{--}, n(1)_- \in Z_+,
$$

$$
1 \le n(2)_{--}  \le n(2)_- < \infty, \ n(2)_{--}, n(2)_- \in Z_+.
$$

 \ We denote analogously by $ \  \pi^+(L) \  $ the set of all rectangles which are {\it \ circumscribed \ }   about
 the set $ \ L: \ \pi^+(L) = \{ L^+  \}, $  where

$$
L^+ = \{ [n(1)^{+}, n(1)^{++}] \otimes [n(2)^{+}, n(2)^{++}]  \}: \ L^+ \supset L, \eqno(4.10)
$$

$$
1 \le n(1)^{+}  \le n(1)^{++} < \infty, \   n(1)^{+},  n(1)^{++} \in Z_+,
$$

$$
1 \le n(2)^{+}  \le n(2)^{++} < \infty, \   n(2)^{+},  n(2)^{++} \in Z_+.
$$

\vspace{4mm}

 \ Let now  $ \  \{ L \} = \{L_{\alpha} \}, \ \alpha \to \alpha_0 \  $ be a sequence or more generally a net of  subsets in $ \  Z_+^2, \  $
and let  $ \  \{ L_- \} = \{L_{-,\alpha} \}, \ \alpha \to \alpha_0 \  $ be arbitrary  correspondent sequence of inscribed rectangles:

$$
L_- = L_-(L_{\alpha}) = \{ [n(1)_{--}, n(1)_-] \otimes [n(2)_{--}, n(2)_- ]  \}: \ L_- \subset L.
$$

 \ Denote

$$
\kappa_- = \kappa_-(L) = \kappa(L,L_-) \stackrel{def}{=}  \frac{|L \setminus L_-|}{|L|^{1/2}}. \eqno(4.11)
$$

 \ We  impose the following two conditions on the sequence (net) of the sets $ \   \left\{ L_{\alpha}, \  L_-(L_{\alpha}) \right\}: $

$$
\min \left[n(1)_{--}, n(2)_{--} \right] \to \infty; \eqno(4.12)
$$

$$
\kappa_- = \kappa_- \left(L_{\alpha} \right) = \kappa_- \left(L_{\alpha}, L_{-,\alpha}\right)  \to 0. \eqno(4.13)
$$

 \vspace{4mm}

 \ {\bf Theorem 4.2.} Suppose that  both the conditions (4.12) and (4.13) are satisfied. Let also  $ \ 0 < W[f](2) <  \infty. \ $ Then
the assertion of theorem 4.1 remains true. \par

 \vspace{4mm}

 \ {\bf Proof.}  It follows from theorem 4.1 that the sequence $ \  S_{L_-} \  $ converges in distribution to the multiple stochastic integral
$ \ S_{\infty}. \ $ Therefore, it  is sufficient to ground the convergence

$$
\forall \epsilon > 0 \ \Rightarrow \ {\bf P} \left( \left|  S_L - S_{L_-}  \right| > \epsilon \right) \to 0, \eqno(4.14)
$$
i.e. estimate the remainder term $ \ S_L - S_{L_-} . \  $\par

 \ Note first of all for this purpose

$$
|L|^{1/2} \cdot S_L = \sum_{k \in L} f(\vec{\xi}) = \sum_{k \in L_-} f(\vec{\xi}) +  \sum_{k \in L \setminus L_-} f(\vec{\xi})
\stackrel{def}{=} \Sigma_1 + \Sigma_2;
$$

$$
S_{L_-} =  \left[\frac{|L|}{|S_{L_-}|}  \right]^{1/2} \cdot S_L \to 1 \cdot  S_{\infty} = S_{\infty},
$$
as long as $ \  |S_{L_-}|/|S_L| \to 1. \  $\par

 \ Further,

$$
\left| \ |L|^{-1/2} \ \Sigma_2 \right|_2 \le W_2[f] \cdot  \frac{|L| - |L_{-}|}{|L|^{1/2}} \le W_2[f] \cdot
\kappa_- \left(L_{\alpha} \right)  \to 0,
$$
 \ Q.E.D. \par

\vspace{4mm}

 \ The case of  circumscribed rectangle, as well as the multidimensional case $ \ d \ge 3 \ $ may be investigated alike the considered one.
Indeed, we  impose as above the following two conditions on the sequence (net) of the sets $ \   \left\{ L^+{\alpha}, \  L^+(L^+{\alpha}) \right\}: $

$$
\min \left[n(1)^+, n(2)^{++} \right] \to \infty; \eqno(4.15)
$$

$$
\kappa^+ = \kappa^+ \left(L_{\alpha} \right) = \kappa^+ \left(\ L^+{\alpha}, L^+(L{^+,\alpha}) \ \right)  \to 0, \eqno(4.16)
$$
where

$$
\kappa^+ = \kappa^+(L) = \kappa(L,L^+) \stackrel{def}{=}  \frac{ \left| L^+ \setminus L \ \right|}{|L|^{1/2}}. \eqno(4.17)
$$

 \vspace{4mm}

 \ {\bf Theorem 4.3.} Suppose that  both the conditions (4.15) and (4.16) are satisfied. Let also  $ \ 0 < W[f](2) <  \infty. \ $ Then
the assertion of theorem 4.1 remains true. \par

\vspace{4mm}

 \section{Multidimensional generalization.}

\vspace{4mm}

 \ Let now $ \  (X_m, \ B_m, \mu_m), \ m = 1,2,\ldots, d, \ d \ge 3  \ $ be a family of probability spaces: $ \  \mu_m(X_m)  = 1; \  $
$ \ X := \otimes_{m=1}^d X_m; \ \xi(m) \ $ be independent random variables having the distribution correspondingly
$ \ \mu_m: \ {\bf P}(\xi(m) \in A_m) = \mu_m(A_m), \ A_m \in B_m; \ $
$  \  \xi_i(m), \ i = 1,2,  \ldots, n(m); \ n(m) = 1,2,  \ldots, \ n(m) < \infty  \  $ be independent  copies of $ \ \xi(m) \ $ and also independent
on the other vectors  $ \  \xi_i(s), s \ne m, \ $ so that all the random variables $ \ \{ \xi_i(m)  \} \ $ are common independent. \par

 \ Another notations, conditions, restrictions and definitions. $ \  L  \subset Z_+^d, \ |L| = \card(L) > 1;  \ j = \vec{j} \in L; \  $

$$
 k = \vec{k} = (k(1), k(2), \ldots, k(d)) \in Z_+^d; \  N(\vec{k}) := \max_{j = 1,2, \ldots,d} k(j); \eqno(5.0)
$$

$ \vec{\xi} := \{\xi(1), \xi(2), \ldots, \xi(n(m)) \}; \   \vec{\xi}_i := \{\xi_i(1), \xi_i(2), \ldots, \xi_i(n(m)) \};  $
$ \ X := \otimes_{i=1}^d X_i, \ f:X \to R \  $ be measurable {\it centered} function, i.e. such that $ \ {\bf E} f(\vec{\xi}) = 0; \ $

$$
S_L[f] := |L|^{-1/2} \sum_{k \in L} f\left(\vec{\xi}_k \right). \eqno(5.1)
$$

 \ The following simple estimate is named as before trivial:

$$
|S_L[f]|_p  \le |L|^{1/2} \ |f|L_p. \eqno (3.0a)
$$

\vspace{4mm}

{\it  Recall that  by-still  hereafter \   } $ \ p \ge 2. \ $ \par

\vspace{4mm}

 \ By definition, as above, the function $ \  f: X \to R  \ $ is said to be degenerate, iff it has the form

$$
f(\vec{x}) = \sum_{\vec{k} \in Z_+^d, \ N(\vec{k}) \le M} \lambda(\vec{k}) \ \prod_{s=1}^d g^{(s)}_{k(s)}(x(s)), \eqno(5.2)
$$
for some integer {\it constant} value $ \ M, \ $  finite or not, where all the functions $ \   g^{(s)}_k(\cdot) \  $ are in turn centered:
$ \  {\bf E} g^{(s)}_k(\xi(k)) = 0. \ $  Denotation: $ \ M = \deg[f]. $ \par

\vspace{4mm}

 \ Define also as in the two-dimensional case for each such a function $ \   f \in D \  $ the following non - negative quasi - norm

$$
||f|| D_p \stackrel{def}{=} \inf \left\{ \ \sum_{\vec{k} \in Z^d_+, \ N(\vec{k}) \le M[f] }
|\lambda(\vec{k})|  \cdot \prod_{s=1}^d  |g^{(s)}_{k(s)}(\xi(s))|_p \  \right\}, \eqno(5.3)
$$
where all the arrays $ \ \{ \lambda( \vec{k}) \} , \ \{ g_j \},  \ $ are taking from the  representation 5.2. \par

 \ The last assertion allows a simple estimate: $ \  ||f||D_p \le || f ||D^o_p, \  $ where

$$
||f|| D^o_p \stackrel{def}{=}  \ \sum_{\vec{k} \in Z^d_+, \ N(\vec{k}) \le M[f] }
|\lambda(\vec{k})|  \cdot \prod_{s=1}^d  |g^{(s)}_{k(s)}(\xi(s))|_p, \eqno(5.3a)
$$
and if we denote

$$
G(p) :=  \prod_{j=1}^d  |g^{(j)}_{k_j}(\xi_j)|_p, \ p \ge 1; \ || \lambda ||_1 := \sum_{\vec{k} \in Z^d_+}  |\lambda(\vec{k})|,
$$
then

$$
||f||D_p \le ||f||D^o_p \le G(p) \cdot ||\lambda||_1. \eqno(5.3b)
$$

\vspace{4mm}

\ Further, let the function  $ \  f \in B_1 \otimes B_2  \otimes \ldots \otimes B_d \ $ be given. The error of a degenerate approximation
of the function $ \ f \ $ by the degenerate ones of the degree $ \ M \ $ will be introduced as before

$$
Q_M[f](B_1 \otimes B_2 \otimes \ldots B_d) \stackrel{def}{=} \inf_{\tilde{f} \in D[M]} ||f - \tilde{f}||B_1 \otimes B_2 \otimes \ldots B_d =
$$

$$
\min_{\tilde{f} \in D[M]} ||f - \tilde{f}||B_1 \otimes B_2 \otimes \ldots \otimes B_d. \eqno(5.4)
$$
  \ Obviously, $ \  \lim Q_M[f] (B_1 \otimes B_2 \otimes \ldots \otimes B_d) = 0, \ M \to \infty. \  $\par

 \ For brevity:

$$
Q_M[f]_p \stackrel{def}{=} Q_M[f](L_p(X_1) \otimes L_p(X_2) \otimes \ldots \otimes L_p(X_d)). \eqno(5.5)
$$

\vspace{4mm}

 \  The function $ \ \tilde{f} \ $ which realized the minimum in (5.4),  not necessary to be unique,
will be  denoted by $ \  Z_M[f](B_1 \otimes B_2 \otimes \ldots \otimes B_d):  \  $

$$
 Z_M[f](B_1 \otimes B_2 \otimes \ldots \otimes B_d):=
\argmin_{\tilde{f} \in D[M]} ||f - \tilde{f}||B_1 \otimes B_2 \otimes \ldots \otimes B_d, \eqno(5.6)
$$
so that

$$
Q_M[f](B_1 \otimes B_2 \otimes \ldots \otimes B_d) = ||f - Z_M[f]||B_1 \otimes B_2 \otimes \ldots \otimes B_d.  \eqno(5.7)
$$

 \ For brevity:

$$
Z_M[f]_p := Z_M [f](L_p(X_1) \otimes L_p(X_2) \otimes \ldots \otimes L_p(X_d) ). \eqno(5.8)
$$

\vspace{4mm}

 \ We deduce analogously to the third section \par

\vspace{4mm}

 \ {\bf Theorem 5.1.}  Let $ \  f = f(x) = f(\vec{x}), \ x \in X  \  $ be arbitrary function from the space
$ \   L_p(X_1) \otimes L_p(X_2) \otimes \ldots \otimes L_p(X_d), \ p \ge 2. \  $ Then

$$
 \sup_{L: |L| \ge 1} |S_L[f]|_p \le W_d[f](p),
$$
 where $ \ W_d[f](p) \stackrel{def}{=} $

$$
 \sup_{L: |L| \ge 1} \inf_{M \ge 1} \left[ \ K^d_R(p) \ || Z_M[f] ||D_p  + |L|^{1/2}  Q_M[f]_p \ \right]. \eqno(5.9)
$$

\vspace{4mm}

 \ {\bf Example 5.1.}  We deduce  alike the example 3.1 as a particular case

$$
 \sup_{ L: |L| \ge 1} \left| S_L[f] \right|_p \le   K_R^d(p) \cdot ||f||D_p,  \eqno(5.9a)
$$
 if of course the right - hand side  of (5.9a) is finite for some value $ \ p, \ p \ge 2. \ $ \par

\vspace{4mm}

 \ As a slight consequence:

$$
\sup_{ L: |L| \ge 1} \left| S_L[f] \right|_p \le  K_R^d(p) \cdot ||f||D^o_p \le K_R^d(p) \cdot G(p) \cdot ||\lambda||_1. \eqno(5.9b)
$$

\vspace{4mm}

 \ {\bf Remark 5.1.} Notice that the last estimates (5.9), (5.9a), and (5.9b) are essentially non-improvable. Indeed, it is
known still in the one-dimensional case $ \ d = 1; \ $ for the multidimensional one it is sufficient to take as a trial
{\it factorizable } function; say, when $ \ d = 2, \ $ one can choose

$$
f_0(x,y) := g_0(x) \ h_0(y), \ x \in X, \ y \in Y.
$$

\vspace{4mm}

\ Let us investigate now the Non-Central Limit Theorem for the multi-indexed sums in the case when $  \   d \ge 3. \  $ We
suppose first of all without loss of generality that in the expression (5.2) the function $ \  g_j(\cdot) \ $ are (common) orthonormal:

$$
{\bf E} g^{(s)}_j(\xi(j)) \ g^{(t)}_l(\xi(l)) = \delta_k^l \cdot \delta_s^t, \ s,t = 1,2,\ldots, d; \eqno(5.10)
$$
$ \ \delta_k^l  \ $ is Kronecker's symbol. \par

 \ Further, we suppose that the set $  \   L, \ L \subset Z_+^d \  $ is a parallelepiped:

$$
L = [1,n(1)] \otimes [1, n(2)] \otimes \ldots \otimes [1, n(d)], \eqno(5.11)
$$
then $ \  |L| = \prod_{m=1}^d  n(m);  \ $  and that

$$
  \min_l n(l) \to \infty.   \eqno(5.12)
$$

\vspace{4mm}

 \ Let us introduce the following array  $ \  \beta = \{  \beta^{(s)}(\vec{k})  \}, \ s = 1,2,\ldots,d; \  \vec{k} \in L \  $ consisting
on   the  (common) independent standard normal (Gaussian) distributed random variables, defined on certain sufficiently rich
probability space. \par

\vspace{4mm}

{\bf Theorem 5.2.} We assert under conditions (5.12) and $ \  W_d[f](2) < \infty, \ $
that the generalized sequence of the r.v. $ \  S_L \ $
converges in distribution as $ \  \min_l(n(l)) \to \infty  \ $ to the multiple stochastic integral  $ \ S_{\infty} \ $ relative the
white Gaussian random measure with non-random square integrable integrand:

$$
S_L \stackrel{d}{\to} S_{\infty} \stackrel{def}{=} \sum_{\vec{k} \in L} \ \lambda(\vec{k}) \cdot \prod_{s=1}^d \beta^{(s)}(k_s), \eqno(5.13)
$$
where the symbol $ \  S_L \stackrel{d}{\to} \nu \ $ denotes the weak (i.e. in distribution) convergence. \par

\vspace{4mm}

 \ {\bf Proof \ }  is at the same as in theorem 4.1 and may be omitted. \par

\vspace{4mm}

 \section{Grand Lebesgue Spaces (GLS).}

\vspace{4mm}

 \ We intend to derive in this section the  uniform relative the  amount of summand $ \ |L| \ $ {\it exponential } bounds
for tail of  distribution of the r.v. $ \  S_L, \ $ based in turn on the moments bound obtained above as well as on the theory
of the so-called Grand Lebesgue Spaces (GLS). We recall now some facts about these spaces and supplement more.  \par

\vspace{4mm}

  \ Let  $ \   (\Omega, \cal{M},  {\bf P})   \ $ be  certain probability space.  Let
also $  \psi = \psi(p), \ p \in [1,b), \ b = \const \in (1, \infty] $ be some bounded
from below: $ \inf \psi(p) > 0 $ continuous inside the semi-open interval  $ p \in [1,b) $
numerical valued function. We can and will suppose  without loss of generality

 $$
\ \inf_{p \in [1,b)} \psi(p) = 1 \eqno(6.0)
$$
and $  \ b = \sup \{ p, \ \psi(p) < \infty  \}, $
 so that   $ \supp \ \psi = [1, b) $ or  $ \supp \ \psi = [1, b]. $ The set of all such a functions will be
denoted by  $ \ \Psi(b) = \{ \psi(\cdot)  \}; \ \Psi := \Psi(\infty).  $ \par

\vspace{4mm}

 \ By definition, the (Banach) Grand Lebesgue Space (GLS)    $  \ G \psi  = G\psi(b) $
consists on all the real (or complex) numerical valued measurable functions
(random variables, r.v.)   $   \  f: X \to R \ $  defined on our probability space and having a finite norm

$$
|| \ f \ || = ||f||G\psi \stackrel{def}{=} \sup_{p \in [1,b)} \left[ \frac{|f|_p}{\psi(p)} \right]. \eqno(6.1)
$$
 \ The function $ \  \psi = \psi(p) \  $ in the definition (6.1) is said to be {\it  the generating function } for this space. \par

\  Furthermore, let now $  \eta = \eta(z), \ z \in S $ be arbitrary family of random variables  defined on any set $ \ z \in S \ $ such that

$$
\exists b = \const\in (1,\infty], \ \forall p \in [1,b)  \ \Rightarrow  \psi_S(p) := \sup_{z \in S} |\eta(z)|_p  < \infty.
$$
 \ The function $  p \to \psi_S(p)  $ is named as a {\it  natural} function for the  family  of random variables $  S.  $  Obviously,

$$
\sup_{z \in S} ||\eta(z)||G\Psi_S = 1.
$$

 \ The family $ \ S \ $ may consists on the unique r.v., say $  \  \Delta: \ $

$$
\psi_{\Delta}(p):= |\Delta|_p,
$$
if of course  the last function is finite for some value $ \  p = p_0 > 1. \  $\par
  \ Note that the last condition is satisfied if for instance the r.v. $ \  \Delta \ $ satisfies the so-called Cramer's
condition; the inverse proposition is not true. \par
 \ The generating $ \ \psi(\cdot) \ $ function in (6.1)  may be introduced for instance as natural one for some family of a functions. \par

 \ These spaces are Banach functional space, are complete, and rearrangement
invariant in the classical sense, see [4], chapters 1, 2; and were investigated in particular in many works, see e.g.  [5], [11]-[12], [25]-[26],
 [28]-[32], [36]  etc. We refer here some used in the sequel facts about these spaces and supplement more. \par

 \ The so-called tail function $ \ T_{f}(y), \ y \ge 0 \ $ for arbitrary (measurable) numerical valued function (random variable, r.v.)
 $ \  f \ $ is defined as usually

$$
T_{f}(y) \stackrel{def}{=}  \max ( {\bf P}(f \ge y), \  {\bf P}(f \le -y) ), \ y \ge 0.
$$

 \ It is known that

$$
|f|^p_p= \int_X |f|^p(\omega) \ {\bf P}(d\omega) = p \int_0^{\infty} y^{p-1} \ T_f(y) \ dy
$$

and if  $  \ f \in G\psi, \ f  \ne 0, $ then

$$
T_{f}(y) \le  \exp \left( -v_{\psi}^*(\ln(y/||f||G\psi)   \right),  \  y \ge e \ ||f ||G\psi, \eqno(6.2)
$$
where

$$
v(p) = v_{\psi}(p) := p \ \ln \psi(p).
$$

 \ Here and in the sequel the operator  (non - linear) $ \ f \to f^* \ $   will denote the famous Young-Fenchel, or Legendre transform

$$
f^*(u) \stackrel{def}{=} \sup_{x \in \Dom(f)} (x \ u - f(x)).
$$

 \ Conversely, the last inequality may be reversed in the following version: if

$$
T_{\zeta}(y)  \le  \exp \left( - v_{\psi}^* (\ln (u/K) \right), \ u \ge e \ K,
$$
and if the auxiliary function  $  \ v(p) = v_{\psi}(p)   $ is positive, finite for all the values $ \ p \in [1, \infty),  $ continuous,
convex and such that

$$
\lim_{p \to \infty} \psi(p) = \infty,
$$
then $  \ \zeta \in G(\psi) \  $ and besides $  \ ||\zeta|| \le C(\psi) \cdot K.  $\par

\vspace{4mm}

 \ Let us consider the so-called {\it exponential} Orlicz space $  \Phi(M) \ $  equipped with an ordinary Luxemburg norm
builded over source probability space with correspondent Young-Orlicz function

$$
M(y) = M[\psi](y) = \exp \left(  v_{\psi}^*(\ln |y|) \ \right), \ |y| \ge e;  \ M(y) = C y^2, \ |y| < e. \eqno(6.2a)
$$
 \ Of course, $ \   C e^2 = \exp \left(v_{\psi}^*(1)  \right). \  $\par

\vspace{4mm}
 \ The  {\it exponentiality} implies in particular that the Orlicz space  $ \ \Phi(M) \ $ is not separable in general case as long as the correspondent
Young - Orlicz function  $  \ M(y) = M[\psi](y) \ $ does not satisfy the $ \ \Delta_2 \ $ condition. \par

 \ The Orlicz-Luxemburg $ \ ||\cdot||\Phi(M) = ||\cdot||\Phi(M[\psi](\cdot)) \  $ and $ \ ||\cdot|| G\psi \ $ norms are quite equivalent:

$$
||f||G\psi  \le C_1 ||f||\Phi(M) \le C_2 ||f||G\psi,
$$

$$
 0 < C_1 = C_1(\psi) < C_2 = C_2(\psi) < \infty. \eqno(6.3)
$$

\vspace{4mm}

  \ {\bf Example 6.0.}  \ Let us consider also the so-called {\it extremal} $  \ \Psi \ - \ $ function $ \ \psi_{(r)}(p), \ $
where $  r = \const \in [1,\infty): $

$$
\psi_{(r)}(p) \stackrel{def}{=} 1,  \ p \in [1,r];
$$
so that the correspondent value $  b = b(r) $  is equal to $  r. $  One can  extrapolate formally this function onto the whole  semi-axis $  R^1_+: $

$$
\psi_{(r)}(p)  := \infty, \ p > r.
$$

 \  The classical Lebesgue-Riesz $ L_r  $ norm for the r.v. $  \eta $ is  quite equal to the GLS norm $  ||\eta|| G\psi_{(r)}: $

$$
|\eta|_r = ||\eta|| G\psi_{(r)}.
$$
 \ Thus, the ordinary Lebesgue-Riesz spaces are particular, more precisely, extremal cases of the Grand-Lebesgue ones. \par

\vspace{4mm}

  \ {\bf Example 6.1.}  For instance, let  $ \psi $  function has a form

$$
\psi(p) = \psi_m(p) = p^{1/m}, \ p \in [1, \infty),  \ m = \const > 0.  \eqno(6.4)
$$

 \ The function $ \ f: X \to R \ $ belongs to the space $ \ G\psi_m: \ $

$$
||f||G\psi_m = \sup_{p \ge 1} \left\{ \  \frac{|f|_p}{p^{1/m}}  \ \right\} < \infty
$$
if and only if the correspondent tail estimate is follow:

$$
\exists  V = V(m) > 0 \ \Rightarrow \   T_f(y)  \le \exp \left\{ -  (y/V(m))^m   \right\}, \ y \ge 0.   \eqno(6.5)
$$

 \  The correspondent  Young-Orlicz function for the space $ \ G\psi_m \ $  has a form

$$
M_m(y) = \exp \left(  |y|^m  \right), \ |y| > 1; \ M_m(y) = e \ y^2, \ |y| \le 1.
$$

 \ There holds for arbitrary  function $ \  f \ $

$$
|| f ||G\psi_m  \asymp ||f||\Phi(M_m) \asymp V(m),
$$
if of course as a capacity of the value $ \ V = V(m) \ $ we understand its minimal positive value from the relation (6.5). \par

 \ The case $ \  m = 2  \ $ correspondent to the so-called subgaussian case, i.e. when

$$
T_f(y) \le \exp \left\{ -  (y/V(2))^2   \right\}, \ y > 0.
$$

\vspace{4mm}

 \ It is presumes as a rule  in addition that the  function $ \  f(\cdot) \ $ has a mean zero:  $ \  \int_X f(x) \ \mu(dx) = 0. \ $
More examples may be found in [6], [25], [26], [28], [29], [30], [31], [36] etc. \par

\vspace{4mm}

 \  We   bring  a more general example, see [26].  Let $  \  m = \const > 1 \  $ and define $  \  q = m' = m/( m-1). \  $
Let also $  \  R = R(y), \ y > 0 \ $ be positive  continuous differentiable {\it  slowly varying  }  at infinity function such that

$$
\lim_{\lambda \to \infty} \frac{R(y/R(y))}{R(y)} = 1. \eqno(6.6)
$$
 \ Introduce a following $ \ \psi \ - \ $ function

$$
\psi_{m,L} (p) \stackrel{def}{=} p^{1/m} R^{-1/(m-1)} \left(  p^{ (m-1)^2/m  }  \right\}, \ p \ge 1, m = \const > 1, \eqno(6.7a)
$$
and  the  correspondent exponential tail function

$$
T^{(m,R)}(y) \stackrel{def}{=} \exp \left\{  - \ y^m \ R^{ m-1} \left(y^{m-1} \right)  \right\}, \ y > 0. \eqno(6.7b)
$$

 \ The following implication holds true:

$$
 0 \ne f \in G\psi_{m,L} \ \Longleftrightarrow \exists C = \const \in (0,\infty), \ T_f(y) \le  T^{(m,R)}(y/C). \eqno(6.8)
$$
 \ A particular cases: $ \  R(y) = \ln^r (y+e), \ r = \const, \ y \ge 0; $ then the correspondent generating  functions has a form

$$
\psi_{m,r}(p) = \ p^{1/m}  \ \ln^{-r}(p), \ p \in [2, \infty),  \eqno(6.9a)
$$
and the correspondent tail function
$$
T^{m,r}(y) = \exp \left\{ \ -  y^m \ (\ln y)^{r}   \ \right\}, \ y \ge e. \eqno(6.9b)
$$
\ For instance, for the Poisson distribution

$$
\psi(p) = \frac{C_R \ p}{ e \ \ln p}, \ p \ge 2,
$$
and

$$
T(y) = \exp (- y \ \ln y), \ y \ge e.
$$

 \ For the exponential distribution $ \  r = 0; \  $ therefore the correspondent $ \ \psi \ - \ $ function  has a form $ \ \psi(p) = p.  \  $
For Gaussian distribution or more generally subgaussian one again $ \ r = 0 \ $ and $ \ m = 2. \ $ \par

 \ More precisely, if

$$
||\xi||G\psi_{m,0} \le 1, \ \Leftrightarrow |\xi|_p \le p^{1/m},  \ p \ge 1,
$$
then

$$
T_{\xi}(y) \le \exp \left\{ -  (me)^{-1} y^m  \right\}, \ y \ge 0.
$$

 \ The inverse conclusion is also true up to multiplicative constant. \par

\vspace{4mm}

 \ {\bf Example 6.2.}  Bounded support of generating function. \par

\vspace{4mm}

 \ Introduce the following tail function

$$
T^{<b,\gamma, R>}(x)  \stackrel{def}{=} x^{-b} \ (\ln x)^{\gamma} \ R(\ln x), \ x \ge e, \eqno(6.10)
$$
where as before $ \   R = R(x), \ x \ge 1 \ $ is positive continuous slowly varying function as $ \ x \to \infty, \ $ and

$$
b = \const \in (1, \infty), \ \gamma = \const > -1.
$$

 \ Introduce also the following (correspondent!) $ \ \Psi(b) \ $ function

$$
\psi^{<b,\gamma, R>}(p) \stackrel{def}{=} C_1(b,\gamma,R) \ (b-p)^{ -(\gamma + 1)/b } \ R^{1/b} \left(  \frac{1}{b-p} \right), \
1 \le p < b. \eqno(6.11)
$$

  \ Let  the measurable function (r.v.) $ \ f(\cdot) \ $ be such that

$$
T_f(y) \le T^{<b,\gamma, R>}(y),  \ y \ge e,
$$
then

$$
|f|_p \le  C_2(b,\gamma, R) \ \psi^{<b,\gamma, R>}(p),  \ p \in [1,b)  \eqno(6.12)
$$
or equivalently

$$
 ||f|| \in  G \psi^{<b,\gamma, R>} \ \Longleftrightarrow \  ||f||G \psi^{<b,\gamma, R>}  < \infty. \eqno(6.13)
$$

 \ Conversely, if the estimate (6.12) holds true, then

$$
T_f(y) \le   C_3(b,\gamma,L)  \ y^{-b} \ (\ln y)^{\gamma + 1} \ R(\ln y), \  y \ge e  \eqno(6.14a)
$$
or equally

$$
T_f(y) \le T^{<b,\gamma + 1, R>}(y/C_4),  \ y \ge C_4 \ e. \eqno(6.14b)
$$

 \  Notice that there is a  logarithmic ``gap''   as $ \  y \to \infty \ $  between the estimations  (6.13) and (6.14).  The
cause of this effect is following: the support on these  $ \  \psi \ - \  $ function is bounded, in contradiction to the other examples. \par

 \ Wherein all the estimates (6.14) and (6.13) are non - improvable, see [26], [28], [32]. \par

\vspace{4mm}

{\bf Example 6.3.}

\vspace{4mm}

 \  Let us consider the following $  \   \psi_{\beta,C}(p) \  $ function

$$
  \psi_{\beta,C}(p)  :=  \exp \left( C p^{\beta} \right), \  C, \ \beta = \const > 0. \eqno(6.15)
$$
 \ Obviously, the r.v. $  \tau $ for which

$$
\forall p \ge 1 \ \Rightarrow \ |\tau|_p \ge \psi_{\beta,C}(p)
$$
does not satisfy  the Cramer's condition. \par

\ Let $ \   \xi \ $ be a r.v. belongs to the $ \  G \psi_{\beta,C}(\cdot) \  $ space:

$$
||\xi||  G \psi_{\beta,C} = 1, \eqno(6.16a)
$$
or equally

$$
|\xi|_p \le \exp \left\{ C p^{\beta} \ \right\}, \  p \in [1, \infty). \eqno(6.16b)
$$

 \  The last restriction is quite equivalent to the following tail estimate

$$
T_{\xi}(y) \le \exp \left(  \ - C_1(C, \beta) \ [  \ln(1 + y)   ]^{1 +1/\beta}  \  \right),  \ y > 0, \eqno(6.17)
$$
and the following Orlicz norm finiteness

$$
||\xi||] \Phi (N_{\beta}) \le C_2(C,\beta) < \infty,
$$
where

$$
N_{\beta}(u) := \exp \left(  \  C_3(C, \beta) \ [  \ln(1 + |u|)   ]^{1 +1/\beta}  \  \right),  \ |u| \ge 1. \eqno(6.18)
$$

\vspace{4mm}

\ {\bf Remark  6.1.} \ These GLS spaces are used, for example, for obtaining of an  {\it exponential estimates} for sums of
independent and dependent random  variables and fields, estimations for non-linear functionals from random fields, theory of Fourier series
and transform, theory of operators  etc., see e.g. [5], [18], [25],  [26], [29], sections 1.6, 2.1-2.5. \par

 \vspace{4mm}

\section{Exponential estimates for tail of multi-index sums.}

 \vspace{4mm}

 \ Suppose in this section that the  conclusions (5.9) and (5.9a)  of theorem 5.1 holds true for all the values $ \ p \ $ belonging to
some non-trivial interval $ \  p \in [1,b),  \  $ finite or not: $ \ b = \const \in (1,\infty]: \  $

$$
 W_d[f](p) < \infty, \ p \in [1,b) \eqno(7.1)
$$
or as a particular case

$$
K_R^d(p) \ ||f||D_p < \infty, \ p \in [1,b). \eqno(7.2)
$$

 \ Denote correspondingly

$$
v_{d,W}[f](p) := p \ \ln \left[ \ W_d[f](p) \ \right]; \ v_{K,d}[f](p) := p \ \ln \left[ K_R^d(p) \ ||f||D_p \ \right], \eqno(7.3)
$$
so that

$$
\sup_{L: |L| \ge 1}  |S_L|_p \le v_{d,W}[f](p), \  \sup_{L: |L| \ge 1}  |S_L|_p \le v_{K,d}[f](p)
$$
or equivalently

$$
\sup_{L: |L|  \ge 1} ||S_L||Gv_{d,W} \le 1, \  \sup_{L: |L|  \ge 1} ||S_L||Gv_{K,d} \le 1. \eqno(7.3a)
$$

 \vspace{4mm}

 \ {\bf Proposition 7.1. } \par

\vspace{4mm}

 \ It follows immediately from the theory of Grand Lebesgue Spaces, namely, from the estimate (6.2), the following exponential tail
estimate

$$
\sup_{L: |L| \ge 1} T_{S_L[f]}(u) \le \exp \left(  - v_{d,W}^*[f](\ln u) \  \right), \ u \ge e, \eqno(7.4)
$$

and following as a particular case

$$
\sup_{L: |L| \ge 1} T_{S_L[f]}(u) \le \exp \left(  - v_{K,d}^*[f](\ln u) \  \right), \ u \ge e. \eqno(7.5)
$$

\vspace {4mm}

 \ The relations (7.4) and (7.5) may be rewritten on the terms of Orlicz spaces, in accordance with (6.2), (6.2a), (6.3), as follows.
Introduce the following correspondent Young-Orlicz functions

$$
N_1(y) = \exp \left(  v_{d,W}^*(\ln |y|), \right), \ |y| \ge e;  \ N_1(y) = C_1 y^2, \ |y| < e;
$$

$$
N_2(y) = \exp \left(  v_{K,d}^*(\ln |y|), \right), \ |y| \ge e;  \ N_2(y) = C_2 y^2, \ |y| < e.
$$
and correspondent Orlicz's norms $ \ ||f||\Phi_{N_j}, \ j = 1,2 \ $ defined for the measurable  functions (random variables) with support  on the source
probability space. \par
 \ The inequalities (7.4) and (7.5) are completely correspondingly equivalent up to multiplicative constants under at the same conditions
to the Orlicz norms estimates

$$
\sup_{L: |L| \ge 1} ||S_L[f]||\Phi(N_1) \le C_3 < \infty,
$$

$$
\sup_{L: |L| \ge 1} ||S_L[f]||\Phi(N_2) \le C_4 < \infty.
$$

\vspace{4mm}

 \ Let us bring some examples. \par

\vspace{4mm}

\ {\bf  Example 7.0. }   Let the function $  \ f: X \to R \ $ be from the representation (5.2):

$$
f(\vec{x}) = \sum_{\vec{k} \in Z_+^d, \ N(\vec{k}) \le M} \lambda(\vec{k}) \ \prod_{j=1}^d g^{(j)}_{k_j}(x_j), \eqno(5.2a)
$$
for some constant integer value $ \ M, \ $  finite or not, where all the functions $ \   g_k^{(j)}(\cdot) \  $ are in turn centered:
$ \  {\bf E} g_k^{(j)}(\xi(k)) = 0. \ $ Recall the  denotation: $ \ M = \deg[f]. $ \par

 \ {\it Suppose here  and in what follows in this section that}

$$
\sum_{\vec{k} \in Z_+^d, \ N(\vec{k}) \le M} | \ \lambda(\vec{k}) \ | \le 1
$$
and that each the (centered) r.v. $ \  g_k^{(j)}(\xi(k)) \  $ belongs to  some $ \  G\psi_k \ - \ $ space  uniformly relative the index $ \ j: \ $

$$
 \sup_j| \ g_k^{(j)}(\xi(k)) \ |_p \le \psi_k(p).
$$

 \ Of course, as a capacity of these functions may be picked the natural functions for the r.v. $ \ g_k(\xi(k)): \  $

$$
\psi_k(p) \stackrel{def}{=}  \sup_j |g_k^{(j)}(\xi(k))|_p,
$$
if the last function is finite for some non-trivial interval $ \ [2, a(k)), \ $ where $ \ a(k) \in (2, \infty]. \ $ \par

 \ Obviously,

$$
|f(\vec{\xi})|_p \le \prod_{k=1}^d \psi_k(p),
$$
and the last inequality  is exact if for instance $ \ M = 1 \ $ and all the functions $ \ \psi_k(p) \ $ are natural  for the  family of the
r.v. $ \ g^{(j)}_k(\xi(k)). \  $\par

\vspace{4mm}

 \ Define the following $ \ \Psi \ -  \ $ function

$$
\zeta(p) = \zeta_d[\vec{\xi}](p) \stackrel{def}{=} K_R^d(p) \cdot \prod_{k=1}^d \psi_k(p).
$$

  \ The assertion of proposition (7.1) gives us the estimations

$$
\sup_{L: |L| \ge 1} ||S_L[f]||G\zeta \le 1 \eqno(7.6a)
$$
and hence

$$
\sup_{L: |L| \ge 1} T_{S_L[f]}(u) \le \exp \left( -v^*_{\zeta}(\ln u) \right), \ u \ge e, \eqno(7.6b)
$$
with correspondent Orlicz norm estimate. \par

\vspace{4mm}

\ {\bf  Example 7.1. } \par

 \ Suppose again that

$$
\sum_{\vec{k} \in Z_+^d, \ N(\vec{k}) \le M} \ | \lambda(\vec{k}) \ | \le 1 \eqno(7.7)
$$
and that the  arbitrary r.v. $ \  g^{(j)}_k(\xi(k)) \  $ belongs uniformly relative the index $ \ j \ $ to the correspondent
$   \  G\psi_{m(k), \gamma(k)} \  $ space:

$$
\sup_j | \ g^{(j)}_k(\xi(k)) \ |_p \le p^{1/m(k)} \ [\ln \ p]^{\gamma(k)}, \ p \ge 2, \ m(k) > 0, \ \gamma(k) \in R, \eqno(7.8a)
$$
or equally

$$
\sup_j T_{g^{(j)}_k(\xi(k))}(u) \le \exp  \left( - C(k) \ u^{m(k)} \ [\ln u]^{  - \gamma(k) }    \right), \ u \ge e. \eqno(7.8b)
$$

 \ Define the following variables:

$$
m_0 := \left[ d + \sum_{k=1}^d 1/m(k)  \right]^{-1}, \ \gamma_0 := \sum_{k=1}^d \gamma(k) - d, \eqno(7.9)
$$

$$
\hat{S}_L = \hat{S}_L[f] := e^d \ C^{-d}_R \ S_L. \eqno(7.10)
$$

 \ We conclude by means of the proposition  7.1

$$
\sup_{ L: |L| \ge 1} \left| \left|  \hat{S}_L \right| \right| G\psi_{m_0, \gamma_0} \le 1 \eqno(7.11a)
$$
and therefore

$$
\sup_{ L: |L| \ge 1} T_{\hat{S}_L}(u) \le \exp \left\{  - C(d,m_0, \gamma_0) \ u^{m_0} \ (\ln u)^{ - \gamma_0}   \right\},
\ u > e, \eqno(7.11b)
$$

or equally

$$
\sup_{L: |L| \ge 1} || \hat{S}_L[f] ||\Phi_{m_0, \gamma_0} \le C_1(d, m_0, \gamma_0) < \infty, \eqno(7.11c)
$$
where the Young-Orlicz function $ \  \Phi_{m,\gamma}(\cdot)  \  $ is defined as follows

$$
\Phi_{m,\gamma}(u)\stackrel{def}{=} \exp \left( \  u^m \ [\ln |u|]^{-\gamma}  \right), \ |u| \ge e.
$$

\vspace{4mm}

{\bf Example 7.2.}

\vspace{4mm}

 \  Let us consider as above the following $  \   \psi_{\beta}(p) \  $ function

$$
  \psi_{\beta,C}(p)  :=  \exp \left( C p^{\beta} \right), \  C, \ \beta = \const > 0, \eqno(7.12)
$$
see example 6.3, (6.15) - (6.17). \par

\ Let $ \  g_k^{(j)}(\xi(k)) \ $ be  centered independent random variables belonging to the certain $ \  G \psi_{\beta,C}(\cdot) \  $ space
uniformly relative the indexes $ \ k,j: $

$$
\sup_j \sup_k ||g_k^{(j)}(\xi(k))||  G \psi_{\beta,C} = 1, \eqno(7.13a)
$$
or equally

$$
\sup_j\sup_k T_{g_k^{(j)}(\xi(k)}(y) \le \exp \left(  \ - C_1(C, \beta) \ [  \ln(1 + y)   ]^{1 +1/\beta}  \  \right),  \ y > 0. \eqno(7.13b)
$$

  \ Then

$$
\sup_{L: \ |L| \ge 1} T_{S_L}(y) \le \exp \left(  \ - C_2(C, \beta) \ [  \ln(1 + y)   ]^{1 +1/\beta}  \  \right),  \ y > 0, \eqno(7.14a)
$$
or equally

$$
\sup_{L: |L| \ge 1} || S_L[f]||  G \psi_{\beta,C_3(C, \beta)} =  C_4(C, \beta) < \infty. \eqno(7.14b)
$$
or equally

$$
\sup_{L: |L| \ge 1} ||S_L[f] ||\Phi(N_{\beta}) \le C_5(C,\beta) < \infty, \eqno(7.14c)
$$
where (we recall)

$$
\Phi_{\beta}(u) = \exp \left(  \  C_6(C, \beta) \ (\ln(1 + |u|))^{1 +1/\beta}  \  \right),  \ |u| \ge 1.
$$

\vspace{4mm}

\ {\bf Example 7.3.}  Suppose now that the each centered random variable $ \  g_k^{(j)}(\xi(k))  \  $  belongs uniformly relative the index $ \ j \ $
to certain $ \  G\psi^{<b(k), \theta(k)>} \  $ space, where $ \  b(k) \in (2, \infty), \ \theta(k) \in R: \ $

$$
\sup_j || \ g_k^{(j)}(\xi(k)) \ || G\psi^{<b(k), \theta(k)>} < \infty,
$$
where  (we recall)

$$
\psi^{<b(k),\theta(k)>}(p) \stackrel{def}{=} C_1(b(k),\theta) \ (b(k)-p)^{ -(\theta(k) + 1)/b(k) }, \ 1 \le p < b(k). \eqno(7.15)
$$

 \ This case is more complicates than considered before. \par

 \ Note that if the r.v. $ \  \eta \  $ satisfies the inequality

$$
T_{\eta}(y) \le C \ y^{-b(k)} \ [\ln y]^{ \theta(k)}, \ y \ge e,
$$
then $ \eta \in G\psi^{<b(k),\theta(k)>}, $ see the example 6.2.\par

 \ One can assume without loss of generality

$$
b(1) \le b(2) \le b(3) \le \ldots b(d). \eqno(7.16)
$$

 \ Denote

$$
\nu_k(p) := \psi^{ < b(k), \theta(k) >}(p), \ b(0):= \min_k b(k), \eqno(7.17)
$$
so that $ \ b(0) = b(1) = $

$$
 b(2) = \ldots = b(k(0)) < b(k(0) + 1) \le \ldots \le b(d), \ 1 \le k(0) \le d; \eqno(7.18)
$$

$$
\Theta :=  \sum_{k=1}^{k(0)} (\theta(k) + 1)/b(0), \eqno(7.18a)
$$

$$
\upsilon(p) = \upsilon_{\vec{\xi}}[f](p) \stackrel{def}{=} \prod_{l=1}^{k(0)} \nu_l(p) =
C \cdot \left[ \ b(0) - p \  \right]^{  - \Theta }, \ 2 \le p < b(0). \eqno(7.19)
$$

 \ Obviously,

$$
K_R^d(p) \ \prod_{k=1}^d \ \nu_k(p) \le C \ \upsilon(p) =  C  \ \left[ \ b(0) - p \  \right]^{  - \Theta },  \
 C = C_d(\vec{\xi}, \vec{b}, \vec{\theta}, k(0)). \eqno(7.20)
$$

 \ Thus, we obtained under formulated above conditions

$$
\sup_{L: |L| \ge 1} |S_L|_p \le C_2 \ (b(0) - p)^{-\Theta}, \ p \in [2, b(0)) \eqno(7.21)
$$
with correspondent tail estimate

$$
\sup_{L: |L| \ge 1} T_{S_L}(y) \le C_3 \ y^{-b(0)} \ [ \  \ln y  \ ]^{ \ b(0) \ \Theta}, \ y \ge e. \eqno(7.22)
$$

\vspace{4mm}

\section{Banach space valued Non-Central Limit Theorem  for multi-index sums.    }

\vspace{4mm}

\ Let $ \  V = \{  v \}  \ $ be certain compact metrizable space; the concrete  choosing  of the distance function will be
clarified below. We consider in this section the case when the ``kernel'' centered function  $  \   f \  $ in (5.1) dependent continuously
also on the additional (deterministic!) parameter $ \  v: \ v \in V: \ f = f(v, \vec{\xi}) \ $ and when the (generalized) sequence
$ \ S_L[f] \ $ converges in distribution as $ \  \min n(k) \to \infty \  $ in the space of continuous functions $ \ C(V) \ $ to the 
non-degenerate {\it parametric} multiple stochastic integral relative the Gaussian white measure, alike the fourth section for the ordinary
one-dimensional axis. \par

 \ The one-dimensional  parametric  case $ \ d = 1 \ $ correspondent to the classical Central Limit Theorem in the space of
continuous functions, i.e. when the sequence $ \ S_L \ $ converges weakly to the Gaussian continuous centered random field,
  [27], [28], [29]. \par

 \  See also the previous articles: [19], [34]. \par
 \ More details. $ \ X := \otimes_{i=1}^d X_i, \ f : V \otimes X \to R \  $ be again measurable (r.v.) for each value $ \ v \in V $
 {\it centered} function, i.e. random process (field) such that $ \ \forall v \in V \ {\bf E} f(v,\vec{\xi}) = 0; \ $

$$
Q_L[f](v) := |L|^{-1/2} \sum_{i \in L} f\left(v, \vec{\xi}_i \right). \eqno(8.1)
$$

 \ In what follows the kernel function $ \ f: V \times X \to R \ $ will be some  degenerate parametric generalization of the the representation (5.2):

$$
f(\vec{x}) = \sum_{\vec{k} \in Z_+^d, \ N(\vec{k}) \le M} \lambda(v,\vec{k}) \ \prod_{j=1}^d g^{(j)}_{k_j}(x_j), \eqno(8.1a)
$$
for some constant integer value $ \ M, \ $  finite or not, where all the functions $ \   g_k^{(j)}(\cdot) \  $ are in turn centered:
$ \  {\bf E} g_k^{(j)}(\xi(k)) = 0. \ $ Retain the  denotation: $ \ M = \deg[f]. $ \par

\ We impose again without loss of generality that in the expression (8.1a) the function $ \  g_j^{(s)}(\cdot) \ $ are (common) orthonormal:

$$
{\bf E} g^{(s)}_j(\xi(j)) \ g^{(t)}_l(\xi(l)) = \delta_k^l \cdot \delta_s^t, \ s,t = 1,2,\ldots, d; \eqno(8.2a)
$$
$ \ \delta_k^l  \ $ is Kronecker's symbol. \par

 \ Further, we suppose that the set $  \   L, \ L \subset Z_+^d \  $ is a parallelepiped:

$$
L = [1,n(1)] \otimes [1, n(2)] \otimes \ldots \otimes [1, n(d)], \eqno(8.2b)
$$
then $ \  |L| = \prod_{m=1}^d  n(m);  \ $  and that

$$
  \min_l n(l) \to \infty.   \eqno(8.2c)
$$

 \ Let us introduce again the following array  $ \  \beta = \{  \beta^{(s)}(\vec{k})  \}, \ s = 1,2,\ldots,d; \  \vec{k} \in L \  $ consisting
on   the  (common) independent  standard normal (Gaussian) distributed random variables, defined on certain sufficiently rich
probability space. Define the following multiple parametric stochastic integral relative the Gaussian random  white measure
with non-random square integrable parametric integrand:

$$
Q_{\infty}(v) \stackrel{def}{=} \sum_{\vec{k} \in L} \ \lambda(v, \vec{k}) \cdot \prod_{s=1}^d \beta^{(s)}(k_s), \eqno(8.2)
$$

 \ We  intend  to ground the weak convergence in the space of continuous functions $ \ C(V) \ $ the random fields $ \  Q_L[f] (\cdot) \ $
to one $ \ Q_{\infty}[f] (\cdot), \  $   write $ \  f = f(v, \vec{\xi}) \in NCLT, \  $
presuming of course its continuity. More detail: for arbitrary  bounded continuous functional $ \ F: C(V) \to R \ $

$$
\lim_{\min_k n(k) \to \infty} {\bf E} F(Q_L[f](\cdot) = {\bf E} F(Q_{\infty}[f](\cdot). \eqno(8.3)
$$

 \ The term (notion) NCLT was introduced for the first time in the article [7]. See also [35]. \par

 \ For the parametric $ \ U \  - \ $ statistics this problem was investigated, e.g., in  [33], [34] by means of martingale representations.
 We will use here some other methods: degenerate approximation approach, metric entropy technique etc. \par

\ We will investigate as a capacity of such a functions $ \ \{f\} \ $ again the degenerate ones.
 By definition, as above, the function $ \  f: V \otimes X \to R  \ $ is said to be {\it parametric degenerate,}  iff it has the form

$$
f(v, \vec{x}) = \sum_{\vec{k} \in Z_+^d, \ N(\vec{k}) \le M} \lambda(v,\vec{k}) \ \prod_{j=1}^d g^{(j)}_{k_j}(x_j), \eqno(8.4)
$$
for some integer {\it constant} value $ \ M, \ $  finite or not, where all the functions $ \   g_k^{(j)}(\cdot) \  $ are in turn centered:
$ \  {\bf E} g^{(j)}_k(\xi(k)) = 0. \ $ \par

 Some new notation:

$$
\sigma = \sigma_{\lambda} :=  \sup_{v \in V} \sum_{\vec{k} \in Z_+^d} |\lambda(v, \vec{k})|; \eqno(8.5)
$$

$$
\rho(v_1, v_2) = \rho_{\lambda}(v_1, v_2) :=   \sum_{\vec{k} \in Z_+^d}  \left| \ \lambda(v_1, \vec{k}) - \lambda(v_2, \vec{k}) \ \right|,
\ v_1, v_2 \in V. \eqno(8.6)
$$

 \ The function $ \ (v_1,v_2) \to \rho(v_1, v_2) \ $ is quasi - distance function on the set $ \ V. \ $ The metric entropy function
$ \  H(V, r, \epsilon)  \  $ for arbitrary such a distance - function $ \  r = r(v_1, v_2)  \ $  is defined as ordinary as a logarithm of a
minimal set of points $ \ \{v(i)\}, \ v(i) \in V, \ $  which  forms an epsilon set in the whole space $ \  (V, r). \ $ Denote also

$$
N(V,r,\epsilon) = \exp H(V,r, \epsilon).
$$
 \  It entails from the Haussdorf's theorem that the value $ \  N(V,r, \epsilon) \  $ is finite for any $ \ \epsilon > 0 \ $ if the space
$ \  (V,r) \ $ is pre-compact set. One can suppose in the sequel that the set $ \ V \ $ relative  the distance $ \ \rho_{\lambda} \ $ will be
complete, therefore will be compact set. \par

\vspace{4mm}

{\bf Theorem 8.1.} (Power level). Suppose that {\it for some certain value}  $  \ p \in [2, \infty) \ \Rightarrow  G(p) < \infty   \  $ and that
the following so-called {\it entropic}  integral converges

$$
I := \int_0^1 N^{1/p} (V,  K_R^d(p) \ G(p) \ \rho_{\lambda}, \epsilon) \ d \epsilon < \infty. \eqno(8.7)
$$
 \ Then the r.v. $ \ f = f(v, \vec{\xi}) \ $  satisfies the  NCLT in the space of continuous functions $ \ C(V, \rho) \ $
and wherein

$$
\sup_{L: |L| \ge 1} {\bf E}  \sup_{v \in V} |Q_L(v)|^p < \infty. \eqno(8.7a)
$$

\vspace{4mm}

 \ {\bf Proof.} The convergence of all the finite - dimensional distributions for the sequence of r.f. $ \  Q_L(v) \  $ to ones  for
$ \  Q_{\infty}(v) \ $ follows immediately from theorem 5.2.  The continuity a.e. of the limiting r.f. $ \ Q_{\infty}(v) \ $ follows from the main
result of an article [35];  continuity almost everywhere  of the r.f. $ \ Q_L(v),  \ $ as well as its weak compactness its
distributions in the space of continuous functions $ \ C(V, \rho) \ $ may be simple deduced from  theorem of Pizier [38], as long as

$$
\sup_L \sup_{v \in V} \left| \  Q_L(v)  \ \right|_p \le K_R^d(p) \  G(p) \ \sigma_{\lambda} < \infty,
$$

$$
\sup_L \left| \  Q_L(v_1) - Q_L(v_2)  \ \right|_p \le K_R^d(p) \ G(p) \ \rho_{\lambda}(v_1, v_2).
$$

\vspace{4mm}

 \ In order to obtain the so-called exponential level for NCLT, we need to introduce some new notations and assumptions. \par

 \ Suppose that each the r.v.  $ \  g^{(s)}_k(\xi(k))  \ $ belongs to some $ \ G\psi \ $ space:

$$
\max_s |g^{(s)}_k(\xi(k))|_p \le \upsilon_k(p), \ 2 \le p < a(k).
$$
Then

$$
G(p) \le \prod_{k=1}^d  \upsilon_k(p), \ 2 \le p < a \stackrel{def}{=} \min_k a(k) \in (2, \infty]. \eqno(8.8)
$$

\ Therefore

$$
\sup_L \sup_{v \in V} \left| \  Q_L(v)  \ \right|_p \le K_R^d(p) \ G(p) \ \sigma_{\lambda}, \eqno(8.9a)
$$

and analogously

$$
\sup_L \left| \  Q_L(v_1) - Q_L(v_2)  \ \right|_p \le K_R^d(p) \ G(p) \cdot \rho_{\lambda}(v_1, v_2). \eqno(8.9b)
$$

 \ Let us introduce a new $ \ G\psi \ - \ $ function $ \ \tau = \tau(p) \ $ as follows

$$
\tau(p) = \tau_{\vec{\xi}}[f](p) :=  K_R^d(p) \ G(p), \ p \in [2, a). \eqno(8.10)
$$
 \ Evidently, $ \  \tau(\cdot) \in \Psi. \  $ Both the estimates (8.9a) and(8.9b) may be rewritten as follows.

$$
\sup_L \sup_{v \in V} \left|\left| \  Q_L(v)  \ \right| \right|G\tau \le \ \sigma_{\lambda}, \eqno(8.10a)
$$

and analogously

$$
\sup_L \left| \left| \  Q_L(v_1) - Q_L(v_2)  \ \right| \right|G\tau \le  \rho_{\lambda}(v_1, v_2). \eqno(8.10b)
$$

 \ Denote

$$
z(y) := \ln \tau(1/y),  \ w(x) := \inf_{y > 0} (xy + z(y)),
$$

$$
J := \int_0^1 \exp \left( w(H(V, \rho, \epsilon))  \right) \ d \epsilon \ - \eqno(8.11)
$$
be again the (generalized) entropic integral. \par

\vspace{4mm}

 \ {\bf Theorem 8.2.} (Exponential level).  Suppose that for some value $ \   a \in (2, \infty] \ \forall p \in [2,a) \ \Rightarrow G(p) < \infty \  $
and that $  \  \sigma(\lambda) < \infty, \ J < \infty. \  $ Then
  the r.f. $ \ f = f(v, \vec{\xi}) \ $  satisfies the  NCLT in the space of continuous functions $ \ C(V, \rho) \ $
and wherein

$$
\sup_{L: |L| \ge 1}  || \sup_{v \in V} Q_L(v)|| G\tau < \infty. \eqno(8.12)
$$

\vspace{4mm}

 \ {\bf Proof } is at the same as above in theorem 8.1; the more general than Pisier's sufficient condition for weak compactness of the
family of the continuous random fields defined on some metrizable compact space
is obtained, e.g., in  the monograph [29], chapter 4, section 4.3. \par

\vspace{4mm}

\ {\bf Remark 8.1.} The condition $ \ J < \infty \ $ is satisfied if for instance the set  $ \  V \  $ is bounded closure of convex open set in the
Euclidean space $ \  R^l \  $ equipped with ordinary norm $ \ |v| \ $ and the distance function $ \ \rho(\cdot, \cdot) \ $  is such that

$$
\rho(v_1, v_2) \le C \ |v_1 - v_2|^{\alpha}, \ \alpha > 0.
$$
\ In this case

$$
N(V,\rho, \epsilon) \le C \ \epsilon^{l/\alpha}, \ \epsilon \in (0,1).
$$

\vspace{4mm}

 \ {\bf Example 8.2.}  Suppose in (8.11) that the function $ \  \tau = \tau(p) \  $ coincides with one $ \  \psi_{m,r}(p), \ m = \const > 0, \ r = \const \in R.  \  $
The condition $  \  J < \infty \ $ in (8.11) takes a form

$$
\int_0^1 H^{1/m} (V, \rho, \epsilon) \  [  \ln \ H(V, \rho, \epsilon) \  ]^{-r} \ d \epsilon  < \infty. \eqno(8.13)
$$

\vspace{4mm}

 \section{Upper bounds for these statistics. }

  \vspace{4mm}

 \  {\bf A.}  A simple lower estimate in the Klesov's (3.4) inequality may has a form

$$
\sup_{L: |L| \ge 1} \left| S^{(2)}_L   \right|_p  \ge \left| S^{(2)}_1  \right|_p =
 \ |g(\xi)|_p \ |h(\eta)|_p, \ p \ge 2, \eqno(9.1)
$$
as long as  the r.v. $ \   g(\xi), \ h(\eta) \ $ are independent. \par

 \ Suppose now that $ \ g(\xi) \in G\psi_1 \ $ and $ \  h(\eta) \in G\psi_2, \  $ where $ \ \psi_j \in \Psi(b), \ b = \const \in (2, \infty]; \ $
for instance $ \ \psi_j, \ j = 1,2 \ $ must be the natural functions  for these  r.v. Put $ \  \nu(p) = \psi_1(p) \ \psi_2(p);  \  $  then

$$
   \nu(p) \le \sup_{L: |L| \ge 1} \left| S^{(2)}_L   \right|_p \le K_R^2(p)  \ \nu(p). \eqno(9.2)
$$

 \ Assume in addition that $ \ b < \infty; \ $ then $ \ K_R^2(p) \le C(b) < \infty. \ $ We get to the following assertion. \par

\vspace{4mm}

{\bf Proposition 9.1.} We deduce under formulated above in this section conditions

$$
 1 \le  \frac{ \sup_{L: |L| \ge 1} \left|S_L \right|_p}{\nu(p)} \le C(b) < \infty, \ p \in [2,b). \eqno(9.3)
$$

\vspace{4mm}

  \ {\bf B. \  Tail approach.}  We will use the example  7.2 (and notations therein. ) Suppose in addition that all the (independent) r.v. $  \ \xi(k) \ $
have the following tail of distribution

$$
T_{\xi(k)}(y) = \exp \left(  \ -  [\ln(1 + y)]^{1 +1/\beta} \ \right), \ y \ge 0, \ \beta = \const > 0,
$$
i.e. an unbounded support. As we knew,

$$
\sup_{L: \ |L| \ge 1} T_{S_L}(y) \le \exp \left(  \ - C_5(\beta,d) \ [  \ln(1 + y)   ]^{1 +1/\beta}  \  \right),  \ y > 0,
$$
see (7.14a).  On the other hand,

$$
\sup_{L: \ |L| \ge 1} T_{S_L}(y) \ge T_{S_1}(y) \ge  \exp \left(  \ - C_6(\beta,d) \  [\ln(1 + y)]^{1 +1/\beta} \ \right), \ y > 0. \eqno(9.4)
$$

\vspace{4mm}

 \ {\bf C. An example.} Suppose as in the example 7.1 that the independent centered r.v. $ \  g^{(j)}_k(\xi(k)) \ $ have the  standard Poisson  compensated
distribution: $ \ \Law(\xi(k) +1 ) = Poisson(1), \ k = 1,2,\ldots,d. \ $ Assume also that in the representation (5.2a)  $ \ M = 1 \ $
(a limiting degenerate case). As long as

$$
|g_k^{(j)}(\xi(k))|_p  \asymp  C \frac{p}{\ln p}, \ p \ge 2,
$$
we conclude by virtue of theorem 5.1

$$
\sup_{L: |L| \ge 1} \left|  \ S_L \right|_p \le C_2^d \ \frac{p^{2d}}{ [\ln p]^{2d}}, \ p \ge 2, \eqno(9.5)
$$
therefore

$$
\sup_{L: |L| \ge 1} T_{S_L}(y) \le \exp \left(  - C_1(d) \ y^{1/(2d)} \ [\ln y]^{2d}  \right), \ y \ge e. \eqno(9.6)
$$

 \ On the other hand,

$$
\sup_{L: |L| \ge 1} \left| S_L  \right|_p  \ge \left|  S_1 \right|_p \ge C_3(d) \ \frac{p^d}{ [\ln p]^d},
$$
and following

$$
\sup_{L: |L| \ge 1} T_{S_L}(y) \ge \exp \left(  - C_4(d) \ y^{1/d} \ [\ln y]^{d}  \right), \ y \ge e. \eqno(9.7)
$$

\vspace{4mm}

 \section{Concluding remarks. }

\vspace{4mm}

\ {\bf A.} \ It is interest by our opinion to generalize obtained in this report results onto the
mixing sequences or onto  martingales, as well as onto the multiple integrals instead sums. \par

\vspace{4mm}

{\bf B.} \ Perhaps, a more general results may be obtained by means of the so-called method of
majorizing measures,  see [1]-[3], [8]-[10], [16], [27], [36], [39]-[43].  \par

\vspace{4mm}

{\bf C. } \ Possible applications: statistics and Monte-Carlo method, alike [13],  [15] etc.\par

\vspace{4mm}

{\bf D.} \ It is interest perhaps  to generalize the assertions of theorems 4.2 and 4.3 onto the sequences of domains $ \ \{ \ L \ \} \ $
tending to ``infinity'' in the van Hove sense, in the spirit of an article [6]. \par

 \vspace{4mm}

\ {\bf E.} A simple qualitative  analysis of limit theorems (proposition 6.4)  show us that the speed of convergence in regular case
is equal to $ \ O \left[\min_s n(s) \right]^{-1/2}. \ $ \par

\vspace{4mm}

 {\bf References.}

 \vspace{4mm}

{\bf 1. Bednorz W.} (2006). {\it A theorem on Majorizing Measures.} Ann. Probab., {\bf  34,}  1771-1781. MR1825156.\par

 \vspace{4mm}

{\bf 2. Bednorz  W.} {\it  The majorizing measure approach to the sample boundedness.} \\
arXiv:1211.3898v1 [math.PR] 16 Nov 2012 \\

\vspace{4mm}

{\bf 3. Bednorz W.} (2010), {\it Majorizing measures on metric spaces.}
C.R. math. Acad. Sci. Paris, (2010), 348, no. 1-2, 75-78, MR2586748. \\

\vspace{4mm}

{\bf 4. Bennet C., Sharpley R. } {\it Interpolation of operators.} Orlando, Academic
Press Inc., (1988). \\

\vspace{4mm}

{\bf 5. Buldygin V.V., Kozachenko Yu.V.} {\it Metric Characterization of Random
Variables and Random Processes. } 1998, Translations of Mathematics Monograph,
AMS, v.188. \\

\vspace{4mm}

{\bf 6.  A.Bulinskii and N.Kwyzhanovskaya.} {\it  Convergence rate in CLT for vector - valued random fields with
self  - normalization.  }  Probability Theory and Mathematical Statistics, Kiev, KSU, 2012, Vol.26 Issue 2, pp. 261-281. \\

\vspace{4mm}

{\bf  7. Dobrushin R.L., Major P. } {\it  Non - Central Limit Theorem  for non-linear  functionals of Gaussian fields.}
Zeitschrift fur Wahrscheinlichkeitstheory und verwandte Gebiete. 1979, V.50, P. 27-52.\\

\vspace{4mm}

{\bf  8. Fernique X.} (1975). {\it Regularite des trajectoires des function aleatiores gaussiennes.}
Ecole de Probablite de Saint-Flour, IV-1974, Lecture Notes in Mathematic. {\bf 480,} \ 1-96, Springer Verlag, Berlin.\\

\vspace{4mm}

{\bf 9. Fernique X.}  {\it Caracterisation de processus de trajectoires majores ou continues.}
Seminaire de Probabilits XII. Lecture Notes in Math. 649, (1978), 691-706, Springer, Berlin.\\

\vspace{4mm}

{\bf 10. Fernique  X.} {\it Regularite de fonctions aleatoires non gaussiennes.}
Ecolee de Ete de Probabilits de Saint-Flour XI-1981. Lecture Notes in Math. 9, {\bf 76,}  (1983), 174, Springer, Berlin.\\

\vspace{4mm}

{\bf 11. A. Fiorenza.} {\it Duality and reflexivity in grand Lebesgue spaces.} Collect.
Math. 51, (2000), 131-148.\\

\vspace{4mm}

{\bf 12.  A. Fiorenza and G.E. Karadzhov.} {\it Grand and small Lebesgue spaces and
their analogs.} Consiglio Nationale Delle Ricerche, Instituto per le Applicazioni
del Calcoto Mauro Picone”, Sezione di Napoli, Rapporto tecnico 272/03, (2005). \\

\vspace{4mm}

{\bf  13. Frolov A.S., Chentzov  N.N.} {\it On the calculation by the Monte-Carlo method definite integrals depending on the parameters.}
Journal of Computational Mathematics and Mathematical Physics, (1962), V. 2, Issue 4, p. 714-718 (in Russian).\\

\vspace{4mm}

{\bf 14.  Gine,E.,  R.Latala,  and  J.Zinn.} (2000). {\it Exponential and moment inequalities for U \ - \ statistics.}
Ann. Probab. 18 No. 4 (1990), 1656-1668. \\

\vspace{4mm}

{\bf 15. Grigorjeva M.L., Ostrovsky E.I.} {\it Calculation of Integrals on discontinuous Functions by means of depending trials method.}
Journal of Computational Mathematics and Mathematical Physics, (1996), V. 36, Issue 12, p. 28-39 (in Russian).\\

\vspace{4mm}

{\bf 16. Heinkel  B.} {\it Measures majorantes et le theoreme de la limite centrale dan le  space C(S).}
Z. Wahrscheinlichkeitstheory. verw. Geb., (1977). {\bf 38}, 339-351.

\vspace{4mm}

{\bf 17. R.Ibragimov and Sh.Sharakhmetov.} {\it  The Exact Constant in the Rosenthal Inequality for Random Variables with Mean Zero.}
 Theory Probab. Appl., 46(1), 127–132. \\

\vspace{4mm}

{\bf 18. T. Iwaniec and C. Sbordone.} {\it On the integrability of the Jacobian under
minimal hypotheses.} Arch. Rat.Mech. Anal., 119, (1992), 129-143.\\

\vspace{4mm}

{\bf  19.  Nazgue Jenish and Ingmar R.Prucha. } {\it  Central Limit Theorems and Uniform Laws of Large Numbers for Arrays of
  Random Fields.  }  PDF Internet publication, 2017.  \\

\vspace{4mm}

{\bf 20. Oleg Klesov.} {\it A  limit theorem for sums of random variables indexed by multidimensional indices. }
Prob. Theory and Related Fields,  1981, {\bf 58,} \ (3),  389-396. \\

\vspace{4mm}

 {\bf 21. Oleg Klesov.} {\it A  limit theorem for multiple sums of identical distributed independent
 random variables. }  Journal of Soviet Mathematics, September 1987, V. 38 Issue 6 pp. 2321-2326.
 Prob. Theory and Related Fields,  1981, {\bf 58, \ (3),}  389-396. \\

\vspace{4mm}

{\bf 22. Oleg Klesov.} {\it   Limit Theorems for Multi-Indexes Sums of Random Variables. } Springer, 2014.\\

\vspace{4mm}

{\bf 23. O. Klesov.} {\it Limit theorems for multi-indexed sums of random variables.}
 Volume 71 of Probability Theory and Stochastic Modeling. Springer Verlag, Heidelberg, 2014. \\

\vspace{4mm}

{\bf 24. Korolyuk V.S., Borovskikh Yu.V.} (1994). {\it Theory of U-Statistics.} Kluwner Verlag, Dordrecht,
(translated from Russian).\\

\vspace{4mm}

{\bf 25. Kozachenko Yu. V., Ostrovsky E.I.}  (1985). {\it The Banach Spaces of
random Variables of subgaussian Type.} Theory of Probab. and Math. Stat. (in
Russian). Kiev, KSU, 32, 43-57. \\

\vspace{4mm}

{\bf 26.  Kozachenko Yu.V.,   Ostrovsky E.,   Sirota L} {\it Relations between exponential tails, moments and
moment generating functions for random variables and vectors.} \\
arXiv:1701.01901v1  [math.FA]  8 Jan 2017 \\

\vspace{4mm}

{\bf 27.  Ledoux M., Talagrand M.} (1991).  {\it Probability in Banach Spaces.} Springer, Berlin, MR 1102015.\\

\vspace{4mm}

{\bf 28. E. Liflyand, E. Ostrovsky and L. Sirota.} {\it Structural properties of Bilateral Grand Lebesgue Spaces. }
Turk. Journal of Math., 34, (2010), 207-219. TUBITAK, doi:10.3906/mat-0812-8 \\

\vspace{4mm}

{\bf 29. Ostrovsky E.I.} (1999). {\it Exponential estimations for Random Fields and
its applications,} (in Russian). Moscow-Obninsk, OINPE.\\

\vspace{4mm}

{\bf 30. Ostrovsky E. and Sirota L.}  {\it Sharp moment estimates for polynomial martingales.}\\
arXiv:1410.0739v1 [math.PR] 3 Oct 2014 \\

\vspace{4mm}

{\bf 31. Ostrovsky E. and Sirota L.}  {\it  Moment Banach spaces: theory and applications. }
HAIT Journal of Science and Engineering C, Volume 4, Issues 1-2, pp. 233-262.\\

\vspace{4mm}

{\bf 32. Ostrovsky E. and Sirota L.} {\it  Schl\''omilch and Bell series for Bessel's functions, with probabilistic applications.} \\
arXiv:0804.0089v1 [math.CV] 1 Apr 2008\\

\vspace{4mm}

{\bf 33. E. Ostrovsky, L.Sirota.} {\it  Sharp moment and exponential tail estimates for U-statistics.  }
arXiv:1602.00175v1  [math.ST]  31 Jan 2016 \\

\vspace{4mm}

{\bf 34. Ostrovsky E. and Sirota L.} {\it Uniform Limit Theorem and Tail Estimates for parametric U-Statistics.}
 arXiv:1608.03310v1 [math.ST] 10 Aug 2016\\

\vspace{4mm}

{\bf 35. Ostrovsky E.I.} {\it  Non-Central Banach space valued limit theorem and applications.} In: Problems of the theory of probabilistic
distributions.  1989, Nauka, Proceedings of the scientific seminars on Steklov's Institute,  Leningrad, V.11, p. 114-119, (in Russian).\\

\vspace{4mm}

{\bf  36. Ostrovsky E. and Sirota L.} {\it Simplification of the majorizing measures method, with development.} \\
arXiv:1302.3202v1  [math.PR]  13 Feb 2013 \\

\vspace{4mm}

{\bf  37. Ostrovsky E. and Sirota L.} {\it  Theory of approximation and continuity of random processes.} \\
arXiv:1303.3029v1  [math.PR]  12 Mar 2013

\vspace{4mm}

{\bf 38.  Pizier G.} {\it  Condition d'entropic assurant la continuite de certains processus et applications a l'analyse harmonique. }
Seminaire d' analyse fonctionnalle, 1980, Exp. 13,  P. 22-25. \\

\vspace{4mm}

{\bf 39. Talagrand  M.} (1996). {\it Majorizing measure: The generic chaining.}
Ann. Probab., {\bf 24,} 1049-1103. MR1825156.\\

\vspace{4mm}

{\bf 40. Talagrand M.} (2005). {\it The Generic Chaining. Upper and Lower Bounds of Stochastic Processes.}
Springer, Berlin. MR2133757.\\

\vspace{4mm}

{\bf 41. Talagrand M.}  (1987). {\it Regularity of Gaussian processes.}  Acta Math. 159 no. 1-2, {\bf 99,} \ 149, MR 0906527.\\

\vspace{4mm}

{\bf 42. Talagrand  M.} (1990).  {\it Sample boundedness of stochastic processes under increment conditions.}  \ Annals of Probability,
{\bf 18,} \ N. 1, 1-49, \ MR10439.\\

\vspace{4mm}

{\bf 43. Talagrand  M.} (1992). {\it A simple proof of the majorizing measure theorem.}
Geom. Funct. Anal. 2, no. 1, 118-125. MR 1143666.\\

\vspace{4mm}

{\bf 44.  M. J. Wichura.} {\it  Inequalities with applications to the weak convergence of
random processes with multi-dimensional time parameters.}
Ann. Math. Statist., 40: 681-687, 1969. \\

\vspace{4mm}

\end{document}